\documentclass[authoryear]{elsarticle}

\usepackage{natbib}
\usepackage{subfig}
\usepackage{amssymb}
\usepackage{amsmath}
\usepackage{color}
\usepackage{url}
\usepackage{booktabs}
\newsavebox{\tempbox}

\journal{Renewable Energy}

\newcommand{\dx}{\,\mathrm{d}x}
\newcommand{\dt}{\,\mathrm{d}t}

\begin{document}

\begin{frontmatter}

\title{Tidal turbine array optimisation using the adjoint approach}

\author[ic,grant]{S.W. Funke\corref{cor1}}
\ead{s.funke09@imperial.ac.uk}
\author[ic,simula]{P.E. Farrell}
\author[ic,grant]{M.D. Piggott}

\address[ic]{Applied Modelling and Computation Group, Department of Earth Science and Engineering, Imperial College London, London, UK}
\address[grant]{Grantham Institute for Climate Change, Imperial College London, London, UK}
\address[simula]{Center for Biomedical Computing, Simula Research Laboratory, Oslo, Norway}

\begin{abstract}
Oceanic tides have the potential to yield a vast amount of renewable energy. Tidal stream generators
are one of the key technologies for extracting and harnessing this potential.  In order to extract an
economically useful amount of power, hundreds of tidal turbines must typically be deployed in an array.
This naturally leads to the question of how these turbines should be configured to extract the maximum
possible power: the positioning and the individual tuning of the turbines could significantly influence the extracted power,
and hence is of major economic interest.  However, manual optimisation is difficult due to legal
site constraints, nonlinear interactions of the turbine wakes, and the cubic dependence of the power on the
flow speed. The novel contribution of this paper is the formulation of this problem as an
optimisation problem constrained by a physical model, which is then solved using an efficient
gradient-based optimisation algorithm.  In each optimisation iteration, a two-dimensional finite element shallow
water model predicts the flow and the performance of the current array configuration.  The gradient of the
power extracted with respect to the turbine positions and their tuning parameters is then computed in a fraction of the time taken for a
flow solution by solving the associated adjoint equations. These equations propagate causality backwards
through the computation, from the power extracted back to the turbine positions and the tuning parameters. This yields the
gradient at a cost almost independent of the number of turbines, which is crucial for any
practical application.  The utility of the approach is demonstrated by optimising turbine arrays in
four idealised scenarios and a more realistic case with up to 256 turbines in the Inner Sound of the
Pentland Firth, Scotland.
\end{abstract}


\end{frontmatter}

\section{Introduction}

With the increasing cost of energy, tidal turbines are becoming a competitive and promising option
for renewable electricity generation. A key advantage of tidal energy is that the power extracted is predictable
in advance, which is highly attractive for grid management. In order to amortise the fixed costs of
installation and grid connection, arrays consisting of hundreds of tidal turbines must typically be deployed
at a particular site.  This raises the question of where to place the turbines within the site and how to tune them individually 
in order to maximise the power output; finding the optimal configuration is of huge importance as it could substantially
change the energy captured and possibly determine whether the project is economically viable.
However, the determination of the optimal configuration is difficult because of the complex flow
interactions between turbines and the fact that the power output depends sensitively on the flow
velocity at the turbine positions. 

This problem has heretofore been addressed in two different ways. One approach is to simplify the
tidal flow model such that the solutions are either available as explicit analytical expressions,
or are extremely fast to compute. This means that the optimum can be analytically derived, or
that the whole parameter space of possible configurations can be rapidly explored.
For example, \citet{bryden2007} and \citet{garrett2008}
optimised simplified models to derive an estimate for the maximum energy that can be extracted from
a tidal basin.  \citet{vennell2010, vennell2012} used simple one-dimensional models to demonstrate the importance of tuning each turbine individually 
to account for the channel geometry, turbine positions, and the tidal forcing. Thus, optimisation of farms is a crucial step needed to achieve
their full potential.
However, \citet{vennell2012c} observes that this optimisation requires many model runs (if performed naively), thus making it computationally
infeasible to use expensive, physically-accurate flow models for this task.
While this approach can provide a coarse estimate for the
power potential of a site, these simplified models cannot accurately capture the complex
nonlinear flow interactions between turbines.

The second approach is to use more complex flow models to accurately predict the tidal flow, the
turbine wakes, and the resulting power output. These models are usually formulated as numerical solutions to partial
differential equations (PDEs). The computational expense of these models prohibits
the exploration of the whole parameter space \citep{thomson2011}.  Consequently, typically only a handful of manually
identified turbine configurations are investigated in a given scenario \citep{adams2011}.  \citet{divett2013} compared the
power output of four different layouts in a rectangular channel by solving the two-dimensional
nonlinear shallow water equations and was able to improve the power outcome by over $50\%$ compared to a
regular layout.  \citet{lee2010} used a three-dimensional model to investigate how the distance
between adjacent rows in a regular array layout impacts the turbine efficiency and showed an
efficiency decay for distances of less than three times the turbine diameter.  While these studies
show the potential of improving the performance by changing the turbine positions, such manual
optimisation guided by intuition and experience becomes difficult in a realistic domain with
complex bottom bathymetry, flow dynamics and hundreds of turbines. 

In this paper, we present a novel technique for maximising the power extraction of array configurations that combines
the physical fidelity of PDE-based flow models with advanced automated optimisation techniques. This approach allows the identification
of optimal solutions in a computationally feasible number of iterations, circumventing the computational limitations noted
in \citet{vennell2012c}. The turbine configuration problem is
formulated as a PDE-constrained optimisation problem, which is a major topic of research in applied mathematics
\citep{gunzburger2003,hinze2009}.  The resulting maximisation problem is solved using a gradient-based optimisation
algorithm that takes orders of magnitude fewer iterations than genetic algorithms or simulated annealing
approaches (see e.g. \citet{bilbao2009}). In this paper, the power extracted by an array configuration is predicted using a two-dimensional
nonlinear shallow water model, which captures the interactions between the geometry, the turbines, and the
flow. The gradient of the power is efficiently computed using the \emph{adjoint technique} of variational
calculus, which solves an auxiliary system that propagates causality backwards through the physical system.
This yields the gradient at a cost almost independent of the number of turbines to be optimised, which is crucial for
the method to be applied to large arrays. This gradient is used by the optimisation algorithm to automatically
reposition the turbines and to adjust their tuning parameters. The flow solution is re-evaluated, and the algorithm iterated until an optimum is
found.

This approach has several key advantages. Firstly, it closes the optimisation loop, by accounting for the
effects of the turbines on the flow field itself. This is necessary to find the actual optimum of the
nonlinear optimisation problem. Secondly, unlike gradient-free methods, the approach requires a relatively small number 
of model evaluations and scales to large numbers of turbines,
which is necessary for the optimisation of industrial arrays. For example, in section~\ref{sec:pentland},
an array of 256 turbines is optimised in a realistic domain at an approximate cost of 200 flow solutions.
Thirdly, the optimisation algorithm can incorporate complex constraints such as minimum separation distances,
bathymetry gradient constraints, and legal site restrictions. Finally, the same mathematical framework
extends naturally to more realistic flow models such as the Reynolds-averaged Navier-Stokes equations, and
to other functionals such as profit or environmental impact.

The approach is implemented in an open-source software framework called OpenTidalFarm; all code
and examples from this paper are available at\newline \mbox{\url{http://opentidalfarm.org}}.

\subsection{Optimisation algorithms}

Optimisation algorithms can be divided into two categories: gradient-free and gradient-based
algorithms.  Gradient-free optimisation algorithms use the functional of interest (in this case,
power extracted by the array) as a black box.
They proceed by evaluating the functional at many points in parameter space and use these values to
decide which areas merit further exploration. While these methods tend to be
robust and can, under certain smoothness conditions, provably find globally optimal solutions
\citep{rudolph1996}, they typically require a very large number of functional evaluations that scales
linearly or superlinearly with the number of parameters to be optimised.  For example,
\citet{bilbao2009} used a genetic algorithm that mimics the process of natural evolution to optimise
the location of $8$ wind turbines.  The algorithm was able to improve the power output by about
$70\%$ compared to the initial layout after $17, 300$ functional evaluations.  This large number of
evaluations clearly introduces a practical upper limit for the number of turbines that can be
optimised.  This difficulty is compounded if a more realistic (and hence more expensive) model is
used.  

By contrast, gradient-based optimisation algorithms use additional information to update the
position in parameter space at each iteration: the first or higher derivatives of the functional of
interest with respect to the parameters.  Depending on the problem, this can lead to a significant
reduction in the number of iterations required compared to gradient-free algorithms, making these
the only feasible choice for large scale optimisation problems \citep{gunzburger2003}.  One caveat
of applying gradient-based optimisation algorithms is that they find only local optima.  This issue
can be circumvented by using hybrid approaches \citep{huang2009}.  The main difficulty of
applying gradient-based methods is that the implementation of the gradient computation can be difficult for complex models,
as it involves differentiating through the solution of a partial differential equation.

One way to obtain the derivative information is to approximate the gradient using finite differences. However,
a major disadvantage of this approach is that a single gradient evaluation requires a large number of
functional evaluations that scales linearly with the number of optimisation parameters. This sets a practical
upper bound on the number of turbines to be optimised, and discards the main advantage of gradient-based
optimisation algorithms. Alternatively, the tangent linearisation of the model (i.e. the derivative of the model evaluated at a particular solution) 
can efficiently compute the derivative of all outputs with respect to a single input, while
the adjoint linearisation can efficiently compute the derivative of a single output with respect to all
inputs~\citep{griewank2008}. For the turbine optimisation problem, we wish to maximise a single output (the power extracted) with
respect to many input parameters (the positions and tuning parameters of the turbines); this means that the adjoint approach is the
natural choice, as the required gradient information can be computed in a number of equation solves that is
independent of the number of turbines. 

The development of adjoint models is generally considered as very complicated
\citep{giles2000,naumann2011}. However, this problem has been solved in recent work 
for the case where the forward model is discretised using finite elements, in the
high-level FEniCS framework \citep{farrell2012}. This allows for the extremely
rapid development of optimally efficient adjoint models, which significantly reduces
the development effort required to implement gradient-based optimisation algorithms for
PDE-constrained optimisation problems \citep{funke2013a}.

To the best of our knowledge, this paper presents the first application of the adjoint method to
the optimisation of turbine arrays. While the examples are
shown in the marine context, it is expected that the presented techniques can also be applied to the
optimisation of wind farms. As the wind turbine layout problem is both closely related and
better studied, we next review techniques proposed for its solution.

\subsection{Wind farm optimisation}\label{sec:wind_farm_optimisation}


Layout optimisation for wind farms has been addressed in numerous studies, most of which are based on
gradient-free optimisation algorithms.  In particular, evolutionary methods~\citep{back1996} are
known to yield good results (\cite{salcedo-sanz2011} and the references therein). These algorithms
mimic the process of natural evolution by considering a population of candidate solutions on which
it executes an evolutionary process to find the ``fittest'' solution. 

A related method is particle swarm optimisation~\citep{kennedy1995}, which considers a
population of candidate solutions called particles, that move through the parameter space and
influence each other to drive the swarm to the best solution.  \citet{wan2010} applied particle
swarm optimisation on an analytical wake model to optimise the location of $39$ wind turbines
and showed that this approach can yield better optimal solutions than genetic algorithms.
Simulated annealing algorithms~\citep{laarhoven1987} are probabilistic optimisation methods that
exploit an analogy between the way in which a metal heats and cools into a minimum energy
crystalline structure (the annealing process) and the search for a minimum in a more abstract system.
\cite{bilbao2009} used an analytical wake model to compare a simulated annealing algorithm with a
genetic optimisation algorithm.  In a scenario with $47$ turbines, the number of model evaluations
was reduced from $1,036,200$ with the genetic algorithm to $61,802$ with simulated annealing.
However, such a large number of evaluations would be infeasible for a more complex PDE-based model.

Few publications solve the layout problem with gradient-based optimisation algorithms.
\cite{lackner2007} optimised the position of two wind turbines by applying a gradient-based 
optimisation algorithm to a simplified energy production model with an analytical expression. 
\cite{huang2009} combined a genetic algorithm with steepest ascent to accelerate convergence to an optimal solution.  With this
additional derivative information, the number of iterations was reduced by approximately an order of
magnitude to less than 300 iterations, for a similar power extraction.  Finally, \citet{fagerfjall2010} showed how mixed integer
linear programming techniques can be used to optimise for both the number and position of turbines.  

All of these publications use very simplified flow models for which the gradient is either available
analytically or can be easily approximated.  Furthermore, none of the reviewed papers use gradient-based
methods with the adjoint technique to find an optimal configuration in a physically realistic
model. 

While this prior research on wind farm optimisation is relevant, there are key differences between
wind and tidal turbine arrays.  Firstly, the flow in a tidal channel is dominantly driven by the predictable
tidal forcing, while wind flow modelling is inherently stochastic and needs to include the temporal
uncertainty in the magnitude and direction of the wind forcing.  Secondly, the ratio of turbine
height to free-surface elevation is significantly different: while the rotor diameter of a wind
turbine is small compared to the height of the atmosphere, tidal turbines typically have diameters
of around $20$~m and are deployed in water depths of approximately $50$~m or less.  This leads to little
undisturbed flow above the turbine which could contribute to wake recovery and thus potentially increases 
the length of the wake compared to wind turbines~\citep{bryden2004,divett2013}. 
Finally, if a tidal farm occupies a large fraction of the channel's width and depth then the presence of turbines significantly affects the flow velocity upstream of the farm~\citep{garrett2005, vennell2010}. 
In wind farms, this is not the case~\citep{vennell2012b}. This interaction adds an additional complication to the optimisation of tidal farms.

The remaining sections are organised as follows. 
In section 2, we formulate the turbine configuration problem as an optimisation problem constrained by the shallow water equations. 
Section 3 discusses the discretisation and implementation, which is then verified in section 4. 
In section 5, we demonstrate the capabilities of the proposed approach on four idealised scenarios with 32 turbines.
Section 6 presents an application of this approach where the positions of up to 256 turbines are optimised in a geometry motivated by the Inner Sound of the Pentland Firth, Scotland.
Finally, we make some concluding remarks in section 7. 

\section{Problem formulation}
In this section we formulate the optimal configuration of turbines in an array as a PDE-constrained optimisation problem in the following abstract form:
\begin{equation}
\begin{aligned}
 \max_{z, m}~& J(z, m) \\
\textrm{subject to }
&  F(z, m) = 0, \\
&    b_l \le m \le b_u, \\
&   g(m) \le 0,
\end{aligned}\label{eq:abstract_optimisation_problem}
\end{equation}
where $ J(z, m)\in \mathbb R$ is the functional of interest, $m$ are the design parameters, $F(z, m)$ is a PDE operator parameterised by $m$ with solution $z$, $b_l$ and $b_u$ are lower and upper bound constraints for the design parameters, and $g(m)$ enforces additional restrictions on the design parameters.

In this work $z = (u, \eta)$ is the solution (horizontal velocity, free-surface displacement) of the shallow water equations written in the form $F(z, m) = 0$, and $m$ contains the configuration (position and/or tuning parameters) of the turbines.
The bounds $b_l$ and $b_u$ are used to enforce that the turbines remain in a prescribed area (here assumed a rectangle for simplicity), while
$g(m)$ is used to enforce a minimum distance between any two turbines (e.g. as a multiple of the turbine diameter).

The optimisation problem~\eqref{eq:abstract_optimisation_problem} can be reduced by using the fact that the constraint equation $F(z, m) = 0$ implicitly maps any choice of $m$ to a unique solution $z$.
Hence, the solution $z$ can be considered as an implicit function of the optimisation parameters: $z \equiv z(m)$ . 
By substituting this operator into the functional of interest, we obtain the reduced optimisation problem: 
\begin{equation}
\begin{aligned}
 \max_{m}~& J(z(m), m) \\
\textrm{subject to }
             & b_l \le m \le b_u, \\
             & g(m) \le 0.
\end{aligned}\label{eq:reduced_abstract_optimisation_problem}
\end{equation}

While it is possible to solve the optimisation problem in unreduced form \eqref{eq:abstract_optimisation_problem}, we choose to solve the reduced form, as
it is usually preferable for time-dependent governing equations \citep{long2012}.

\subsection{The design parameters}
For the turbine optimisation problem considered here, the design parameter $m$ is a vector containing the positions, and optionally the tuning parameters, of the turbines.

If the turbine tuning parameters are fixed, then the design parameters contain only the $(x, y)$ positions of the $N$ turbines encoded in the form:
\begin{equation*}
\begin{tabular}{ l | c | c | c | c | c | c | c |}
\cline{2-8}
  $m =$ & $x_1$ & $y_1$ & $x_2$ & $y_2$ & \dots & $x_N$ & $y_N$ \\
\cline{2-8}
\end{tabular}.
\end{equation*}
If additionally the turbines are to be individually tuned, then the vector $m$ is extended to contain the friction coefficient $K_i$ of each turbine:
\begin{equation*}
\begin{tabular}{ l | c | c | c | c | c | c | c | c | c | c | c|}
\cline{2-12}
$m =$ & $K_1$ & $K_2$ & \dots & $K_N$ & $x_1$ & $y_1$ & $x_2$ & $y_2$ & \dots & $x_N$ & $y_N$ \\
\cline{2-12}
\end{tabular}.
\end{equation*}
This could be further generalised to account for any number of turbine parameters.

\subsection{The PDE constraint} 
The constraint equation $F(z, m) = 0$ enforces the laws of physics in the optimisation problem \eqref{eq:abstract_optimisation_problem}. 
For gradient-based optimisation, the constraint equation must fulfill some properties. 
Firstly, for every $m$ it must yield a unique solution $z$ so that $z$ can be written as an implicit function of $m$.
Secondly, $F$ must be differentiable and its derivative with respect to $z$ must be continuously invertible. 

In this work, the physical laws are modelled by the nonlinear shallow water equations. 
Let $\Omega$ be the domain of interest and let $(0, T)$ be the simulation period. Then the equations read:
\begin{equation}
\begin{aligned}
  \kappa \frac{\partial  u}{\partial t} +  u \cdot \nabla  u - \nu \nabla^2  u  + g \nabla \eta + \frac{c_b + c_t(m)}{H} \| u\|  u & = 0, \\
\kappa \frac{\partial \eta}{\partial t} + \nabla \cdot \left(H u\right) & = 0, 
\label{eq:shallow_water_equations}
\end{aligned}
\end{equation}
where the unknowns $u:\Omega
\times (0, T] \to \mathbb R^2$ and $\eta:\Omega \times (0, T] \to \mathbb R$ are the depth-averaged
velocity and the free-surface displacement, respectively, $H:\Omega \to \mathbb R$ is
the water depth at rest, $g\in \mathbb R$ is the acceleration due to gravity, $\nu \in \mathbb R$ is the
viscosity coefficient, and $c_b: \Omega \to \mathbb R$ and $c_t(m): \Omega \to \mathbb R$ represent the quadratic
bottom friction and the turbine parameterisation, respectively.  The parameter $\kappa \in \{0, 1\}$
specifies if the stationary ($\kappa = 0$) or the non-stationary problem ($\kappa = 1$) is
considered; in the stationary case the time-dependency of the variable definitions above can be
neglected.

The boundary conditions are as follows: 
on the domain inflow boundary $\partial \Omega_{\textrm{in}}$ a Dirichlet boundary condition is
applied to the velocity.
On the outflow boundary $\partial \Omega_{\textrm{out}}$ the free-surface displacement is set to zero.
Elsewhere, a no-slip or a no-normal flow boundary condition with a free-slip condition for the tangential components is imposed.
Note that these boundary conditions imply that the effect of the array on the free-stream flow is negligible,
which is not the case for large farms \citep{vennell2010}: for large arrays, it would be more appropriate to force the model
with tidal boundary conditions on the free surface displacement or to ensure that the size of the domain includes the far-field region around the farm.

\subsection{The turbine parameterisation} \label{sec:turbine_parameterisation} 

\begin{figure}
    \centering
    \includegraphics[width=0.5\textwidth]{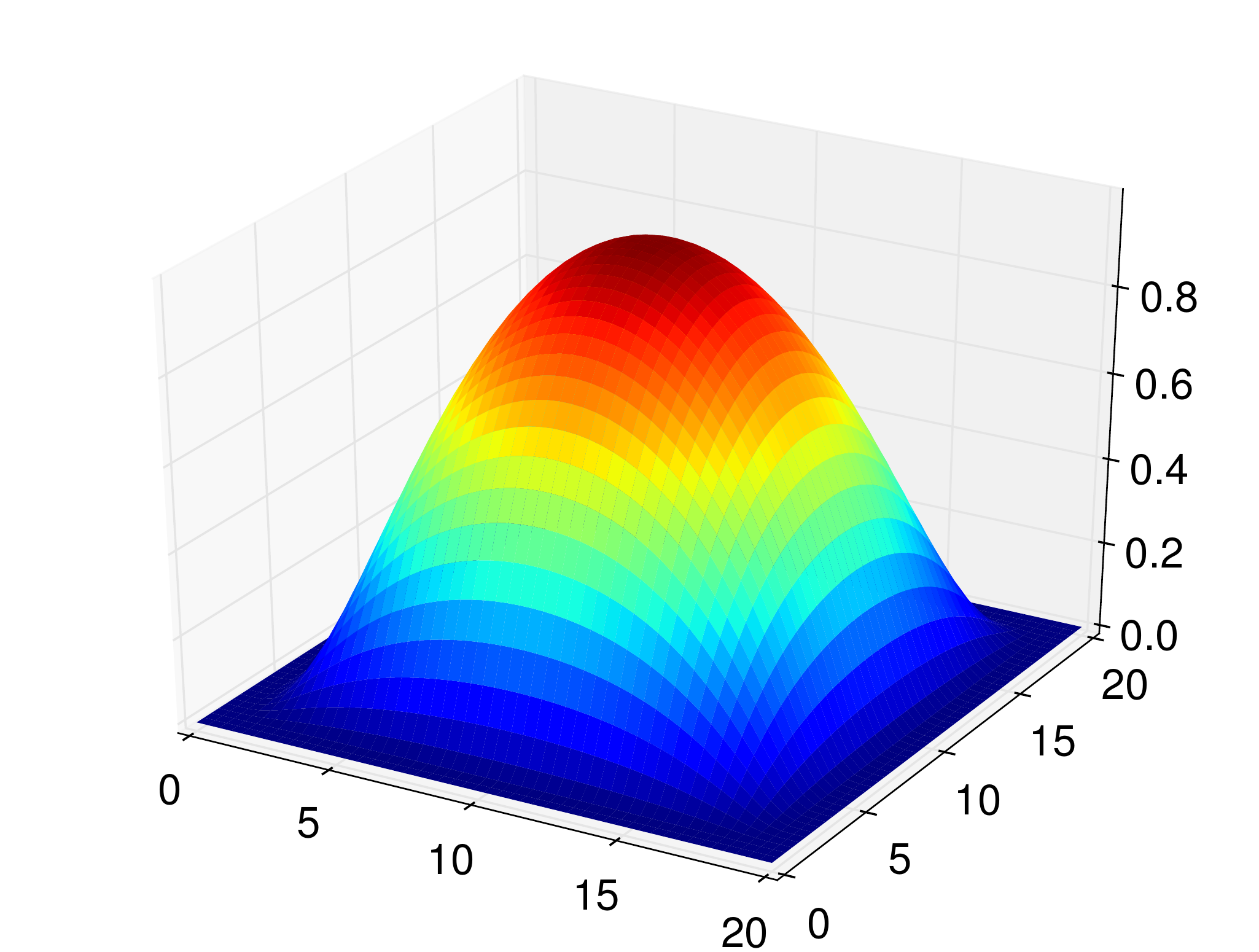}
    \caption{The turbine is parameterised by a smoothly increasing friction coefficient towards the turbine centre, given by the bump function~\eqref{eq:1Dbump_function} multiplied in the $x$ and $y$-direction.}\label{fig:turbine_plot}
\end{figure}

A turbine is modelled here via an increased bottom friction over a small area representative of an individual turbine \citep{divett2013}. 
The sum of the bottom friction associated with all turbines is denoted as $c_t(m)$ in equation~\eqref{eq:shallow_water_equations}.
This individual turbine approach is in contrast to previous work (e.g. \citet{garrett2005,draper2009}), 
where uniformly increased bottom drag over the array area was used to parameterise groups of turbines. 
In \citet{divett2013}, the authors set the friction to a constant value at the turbine locations and zero everywhere else.
This constant parameterisation is problematic in the context of gradient-based optimisation because the friction becomes a non-differentiable function of the turbine position.
For this reason, the turbine parameterisation used in this work smoothly increases the friction value at each turbine position.
The associated friction function is constructed from bump functions, i.e. smooth functions with finite support.
A bump function in one dimension is:
\begin{equation}
\psi_{p, r}(x) \equiv
\begin{cases}
  e^{1-1/(1-\|\frac{x - p}{r}\|^2)} &  \mbox{ for } \|\frac{x - p}{r}\| < 1, \\
  0 & \textrm{ otherwise,}
\end{cases}\label{eq:1Dbump_function}
\end{equation}
where the two parameters $p$ and $r$ are the centre and the support radius of the bump function, respectively.
A two-dimensional bump function is obtained by multiplying equation~\eqref{eq:1Dbump_function} in both independent dimensions.
With that, the friction function of a single turbine parameterised by a friction coefficient $K_i$
centred at a point $(x_i, y_i)$ is given by:
\begin{equation}
 C_{i}(m)(x, y) \equiv K_i \psi_{x_i, r}(x) \psi_{y_i, r}(y).
\label{eq:single_turbine_coefficient}
\end{equation}
A plot of the resulting friction for $K_i=1$, $x_i = 10$, $y_i = 10$ and $r=10$ is shown in figure \ref{fig:turbine_plot}.

The turbine friction function $c_t$ in the governing PDE \eqref{eq:shallow_water_equations} is defined to be the sum of the friction functions~\eqref{eq:single_turbine_coefficient} for all $N$ turbines:
\begin{equation}
 c_t(m) \equiv \sum_{i=1}^{N} C_{i}(m),
\label{eq:turbine_function}
\end{equation}
where a single value for $r$ is used based on the assumption that the deployed turbines are of equal
size.  

Note that turbine properties can be calibrated by modifying the friction parameter $K$ in
equation~\eqref{eq:single_turbine_coefficient}.  For example, the amount of energy that is extracted
from an individual turbine (e.g. due to different pitch settings of the turbine blades) can be
controlled. As a further extension this could be used to handle cut in/out velocities in which the turbines are
operational, but this is not considered in this work.
Finally, other more sophisticated turbine parameterisations such as extensions of actuator disc theory could be employed \citep{roc2013}.

\subsection{The functional of interest} \label{sec:functional}

The functional of interest ${J}$ in equation~\eqref{eq:abstract_optimisation_problem} defines the value of interest that is to be maximised.
For gradient-based optimisation, ${J}$ must be differentiable.

A natural choice for the functional of interest is the time-averaged power extracted due to the increased friction by the turbines~\citep{vennell2012b, sutherland2007}.
In the non-stationary case ($\kappa = 1$) this is expressed as:
\begin{equation}
 J(u, m) = \frac{1}{T} \int_0^T \hspace{-1.75mm} \int_{\Omega} \rho c_t(m) \|  u \|^3 \dx \dt, \label{eq:functional_time_dependent}
\end{equation}
where $\rho$ is the fluid density. 
In the stationary case ($\kappa = 0$) the functional is defined to be the power extracted from the increased friction:
\begin{equation}
  J(u, m) = \int_{\Omega} \rho c_t(m) \|  u \|^3 \dx. \label{eq:functional_steady}
\end{equation}
Note that this value represents kinetic power extraction rather than electrical power generation, since it does not incorporate losses due to the turbine support structures and the conversion to electricity.

More advanced functional choices could instead maximise the profit of the turbine farm, by including
installation and service costs depending on turbine size and the deployment location \citep{adams2011}. Another
alternative would be to incorporate potential environmental impacts.  These are not considered in this
study.

\subsection{Box and inequality constraints}\label{sw_opt:bounds_and_constraints}
The box and inequality constraints in the generic optimisation problem~\eqref{eq:abstract_optimisation_problem} are used to define the feasible values for the optimisation variables. 
In the context of the turbine layout problem, a typical condition is to restrict the area in which the turbines may be placed to the development site.
The numerical examples in this work have rectangular shaped deployment sites and therefore box constraints are sufficient to enforce this restriction.  
For more general site shapes appropriate inequality constraints can be used instead.

Another common condition is to ensure that individual turbines do not overlap. 
This is implemented by enforcing a minimum distance $d_{\textrm{min}}$ between any two turbines:
\begin{equation}
 \|p_i - p_j\|^2 \ge d_{\textrm{min}}^2 \qquad \forall~i, j: 1 \le i < j \le N. \label{eq:minimum_distance_constraint}
\end{equation}

In order for the optimisation to be well-posed, the constraints must satisfy a constraint qualification \citep{hinze2009}. The box constraints
and inequality constraint \eqref{eq:minimum_distance_constraint} are concave functions, and hence satisfy the Concave Constraint
Qualification (CCQ) \citep[theorem 5.4]{carter2001}.

More advanced constraints could for example enforce a minimum or maximum deployment depth or limit the maximum
local bathymetry steepness where turbines may be installed, but these are not investigated in this work.

\section{Numerical setup}\label{sec:turbine_optimisation_discretisation}
\subsection{Optimisation algorithm}
A typical gradient-based optimisation algorithm implements the following iteration:
\begin{itemize}
 \item Choose an initial guess $m^0$ for the design parameters. 
 \item for $i = 0, 1, ...$
  \begin{enumerate}
   \item Solve the forward problem $F(z^i, m^i) = 0$ for $z^i$ and evaluate the functional of interest $J(z^i, m^i)$.
   \item Compute the functional gradient $\mathrm{d}J/\mathrm{d}m(z^i, m^i)$.
   \item Stop if the termination criteria are fulfilled.
   \item Find improved design parameters $m^{i+1}$ using the results of 1 and 2.
  \end{enumerate}
\end{itemize}
Optimisation algorithms differ mainly in their implementation of step 4. The main task is to
identify an improved parameter choice so that the algorithm converges quickly to the optimal
solution while satisfying the imposed constraints.

In this paper we solve the turbine configuration problem using sequential quadratic programming (SQP),
which is considered to be one of the most efficient optimisation algorithms \citep{boggs1996}.  The
implementation used here is the SLSQP algorithm available through the SciPy optimisation package~\citep{scipy} and is described in
detail in \citet{kraft1988}. The optimisation problem was formulated and solved with the PDE-constrained optimisation framework described in
\citet{funke2013a}.


The SQP implementation used is not scale-invariant, i.e. scaling the functional of interest can
impact the convergence of the algorithm (even though it does not change the optimal configuration).
Preliminary numerical investigation found that such rescaling was necessary to achieve fast convergence. 
Therefore, for each numerical experiment presented here the problem was internally rescaled such that the maximum absolute value of the initial functional gradient with respect to the turbine positions was ten times the turbine radius\footnote{More precisely, we compare the \emph{Riesz representer} of the functional gradient with respect to turbine positions, which has the same units of length as the turbine radius.}.

\subsection{The functional gradient computation with the adjoint approach}\label{sec:adjoint_equation}
The second step of the optimisation algorithm requires the computation of the functional derivative with respect to the optimisation parameters $\mathrm{d}J/\mathrm{d}m$.
Its efficient computation is achieved using the adjoint approach. For reasons of brevity, only the key equations are
mentioned; for fuller explanations, see \citet{gunzburger2003, giles2000, farrell2012}.

Let $|J|$ be the number of functional outputs (in this case 1, the power), let $|m|$ be the number of parameters, and let $|z|$ be the number of degrees of freedom of the discretised state vector.
The adjoint approach computes the gradient in two steps. 
Firstly, the \emph{adjoint equation} is solved to obtain the adjoint solution $\lambda$ (with the matrix dimensions indicated below the braces): 
\begin{equation*}
\underbrace{\frac{\partial F}{\partial z}^*}_{|z| \times |z|}
\underbrace{\lambda\vphantom{\frac{\partial F}{\partial z}^*}}_{|z| \times |J|}
 = 
\underbrace{\frac{\partial J}{\partial  z}^*}_{|z| \times |J|},
\end{equation*}
where $*$ denotes the Hermitian transpose. The adjoint solution plays the role of the Lagrange multiplier
enforcing the PDE constraint in the Lagrangian associated with the optimisation problem \eqref{eq:abstract_optimisation_problem}.
Secondly, the \emph{functional gradient} is obtained by:
\begin{equation}
\underbrace{\frac{\mathrm{d}J}{\mathrm{d}m}}_{|J| \times |m|} = 
\underbrace{-\lambda^*\vphantom{\frac{\mathrm{d}J}{\mathrm{d}m}}}_{|J| \times |z|}
\underbrace{\frac{\partial F}{\partial m}}_{|z| \times |m|}
+ \underbrace{\frac{\partial J}{\partial m}}_{|J| \times |m|}. 
  \label{eq:grad_computation_eq}
\end{equation} 
The functional gradient is the key piece of information that makes the optimisation of large numbers of turbines tractable.

For the optimisation problem considered here, the adjoint shallow water equations are:
\begin{equation}
\begin{aligned}
- \kappa \frac{\partial  \lambda_u}{\partial t} + (\nabla  u)^*  \lambda_u - \left(\nabla \cdot  u \right)  \lambda_u -  u \cdot \nabla  \lambda_u - \nu \nabla^2  \lambda_u \\
- H \nabla \lambda_\eta + \frac{c_b + c_t(m)}{H} \left(\| u\|  \lambda_u + \frac{ u \cdot  \lambda_u}{\| u\|}  u\right) & = \frac{\partial J}{\partial  u}^*, \\
- \kappa \frac{\partial \lambda_\eta}{\partial t} - g \nabla \cdot  \lambda_u & = 0. \\
\end{aligned}\label{eq:turbine_optimisation_continuous_adjoint_equations}%
\end{equation}
where $\lambda \equiv (\lambda_{u}, \lambda_\eta)$ is the vector containing the unknown adjoint
velocity and adjoint free-surface displacement, respectively.  The derivation of these equations can
be found in~\citet[appendix C]{funke2012}.  The non-stationary adjoint equations ($\kappa = 1$) have
a final-time condition for the adjoint velocity and free-surface displacement; this condition is
homogeneous, as the functional $J$ has no term evaluated at the end of time.  The adjoint equations
are solved from the final time to the initial time, propagating information backwards in time.  The
boundary conditions for the adjoint equations are the homogeneous versions of the boundary conditions
of the forward equations, again because the functional has no boundary integral terms. The
functional derivative $({\partial J}/{\partial  u})^*$ appears as the source term for the adjoint velocity equation and 
is easy to evaluate as the functional is available as an analytical expression. 
Note that the adjoint equations are linear while the forward equations are nonlinear, and therefore solving the adjoint equations is
typically much cheaper. If the time-dependent equations are solved, the entire forward trajectory is required to assemble
the adjoint equations; if the forward trajectory is too large to store at once, then a checkpointing algorithm must be used
\citep{griewank2000,farrell2012}.

In this work, rather than deriving, discretising and implementing the adjoint equation~\eqref{eq:turbine_optimisation_continuous_adjoint_equations} by hand, we apply the
high-level algorithmic differentiation approach described in \citet{farrell2012}. This efficiently and automatically derives and
implements the discrete adjoint model from the implementation of the forward model~\eqref{eq:to_spatially_discretised_shallow_water_eq}, without user intervention. This significantly
reduces the effort and expertise required to implement adjoints of complex nonlinear forward models.

The second step of the adjoint approach (equation~\eqref{eq:grad_computation_eq}) evaluates the functional gradient using the adjoint
solution $\lambda$.
This step only requires the computation of a matrix-vector product and consequently its computational cost is
negligible. This allows for the computation of the gradient of the functional with only the solution of one adjoint
system. While the cost of the matrix-vector product technically depends on the number of turbines, in practice the cost
of the adjoint PDE solution dominates, and so the adjoint technique yields the gradient at a cost independent of the
number of turbines. This is a key property if many turbines are to be optimised.

\subsection{Discretisation}
The governing PDEs~\eqref{eq:shallow_water_equations} are discretised with the finite element method. 
The weak form is derived by multiplying the equations with test functions $(\Psi, \Phi)$ from
suitable function spaces, integrating over the computational domain and applying integration by
parts to selected terms.
The weak form of equations~\eqref{eq:shallow_water_equations} is:
find $(u, \eta)$ such that $\forall~(\Psi, \Phi)$:
\begin{equation}
\begin{split}
  \kappa \left\langle \frac{\partial  u}{\partial t},  \Psi \right\rangle_\Omega + \left\langle  u \cdot \nabla  u,  \Psi \right\rangle_\Omega + \nu \left\langle\nabla  u, \nabla  \Psi\right\rangle_\Omega & \\
 + g \left\langle \nabla \eta,  \Psi \right\rangle_\Omega  + \left\langle \frac{c_b + c_t(m)}{H} \| u\|  u,  \Psi
 \right\rangle_\Omega &=  ~0, \\
\kappa \left\langle \frac{\partial \eta}{\partial t}, \Phi \right\rangle_\Omega - \left\langle H u, \nabla \Phi \right\rangle_\Omega + \left\langle H  u \cdot  n , \Phi\right\rangle_{\partial \Omega_{\textrm{in}} \cup \partial \Omega_{\textrm{out}}} &=  ~0,
\end{split}\label{eq:to_spatially_discretised_shallow_water_eq}
\end{equation}
where $\left\langle \cdot, \cdot \right\rangle_{\Omega}$ denotes the $L^2(\Omega)$ inner product.
The no-normal flow/free-slip conditions on $\partial \Omega \setminus \left(\partial \Omega_{\textrm{in}} \cup \partial \Omega_{\textrm{out}}\right)$ are weakly imposed by excluding the associated surface integrals.
The Dirichlet boundary conditions are strongly imposed by restricting the function spaces to functions that yield the correct boundary values.

The discretised problem is obtained by choosing discrete function spaces in the weak
formulation~\eqref{eq:to_spatially_discretised_shallow_water_eq}.  In this work, these are
constructed from a suitable triangulation of the computational domain $\Omega$ using the Taylor-Hood
finite element pair which uses piecewise quadratic functions on each triangle for the velocity and
piecewise linear functions on each triangle for the free-surface displacement \citep{hood1973}. 

If the non-stationary problem is considered ($\kappa =1$), then the spatially discretised
equations~\eqref{eq:to_spatially_discretised_shallow_water_eq} must also be discretised in time.  In
this paper, the time discretisation was performed using the implicit Euler method, due to its
unconditional stability and simplicity.

Finally, the integral evaluation of the functionals of interest~\eqref{eq:functional_time_dependent} and \eqref{eq:functional_steady} 
used the same quadrature rules as used in the underlying finite element formulation of the problem.

\section{Verification}

\subsection{Verification of the forward model}
The shallow water model implementation was verified by order-of-convergence analysis.  The
analytical solution is constructed using the method of manufactured solutions~\citep{salari2000,
roache2002}: a desired analytical solution is chosen and then substituted into the governing PDE,
which yields a non-zero remainder.  By adding this remainder as a source term to the governing
equations, the selected solution becomes an analytical solution to the modified PDE. 

For the following tests, the analytical solution consists of sinusoidal functions for both the
velocity and the free-surface displacement:
\begin{equation}
\begin{aligned}
 u_{\mathrm{exact}}(x, y, t) & = 
\begin{pmatrix}
\eta_0\sqrt{gH^{-1}}\cos\left(kx-\sqrt{gH}kt\right) \\
0
\end{pmatrix}, \\
\eta_{\mathrm{exact}}(x,y,t) & = \eta_0\cos\left(kx-\sqrt{gH}kt\right),
\end{aligned}\label{eq:to_mms_test}%
\end{equation}
with $k = \pi /640$~m$^{-1}$, $\eta_0 = 2$~m, $H = 50$~m, $\nu = 3\mathrm{~m}^2\mathrm{s}^{-1}$, $c_t =
0$, $c_b = 0.0025$,  $g = 9.81\mathrm{~ms}^{-2}$ and a final time of $T = \pi/(\sqrt{gH}k)
\approx 28.9$~s which corresponds to half a wave cycle.  The computational domain $\Omega$ is
defined to be a rectangle of size 640 m $\times$ 320 m. 
Following the method of manufactured solutions, the
functions~\eqref{eq:to_mms_test} are substituted into the shallow water equations~\eqref{eq:shallow_water_equations}
and the non-zero remainders are added as source terms. This ensures that \eqref{eq:to_mms_test} is a solution of this modified system.

\begin{figure}[t]
  \centering
     \subfloat[The spatial order of convergence]
     {
        \centering
          \includegraphics[width=0.49\textwidth]{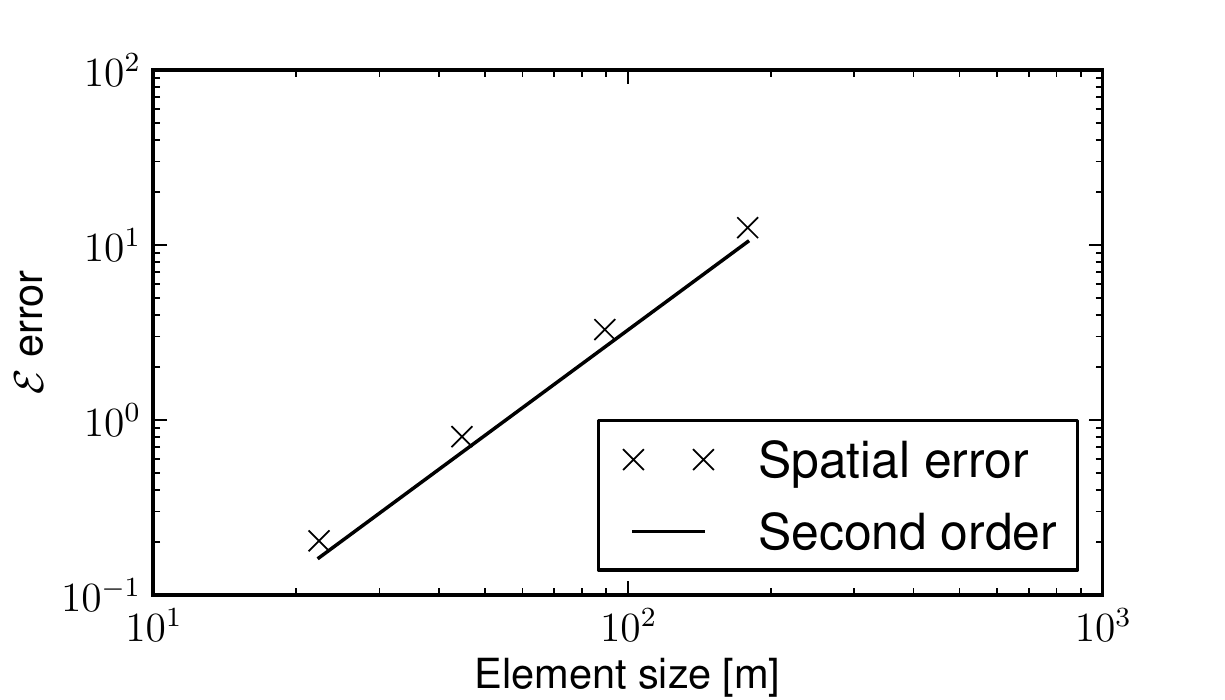}
          \label{fig:results/mms_scenario/spatial_rate_of_convergence} 
         }
         \subfloat[The temporal order of convergence]{
        \centering
           \includegraphics[width=0.49\textwidth]{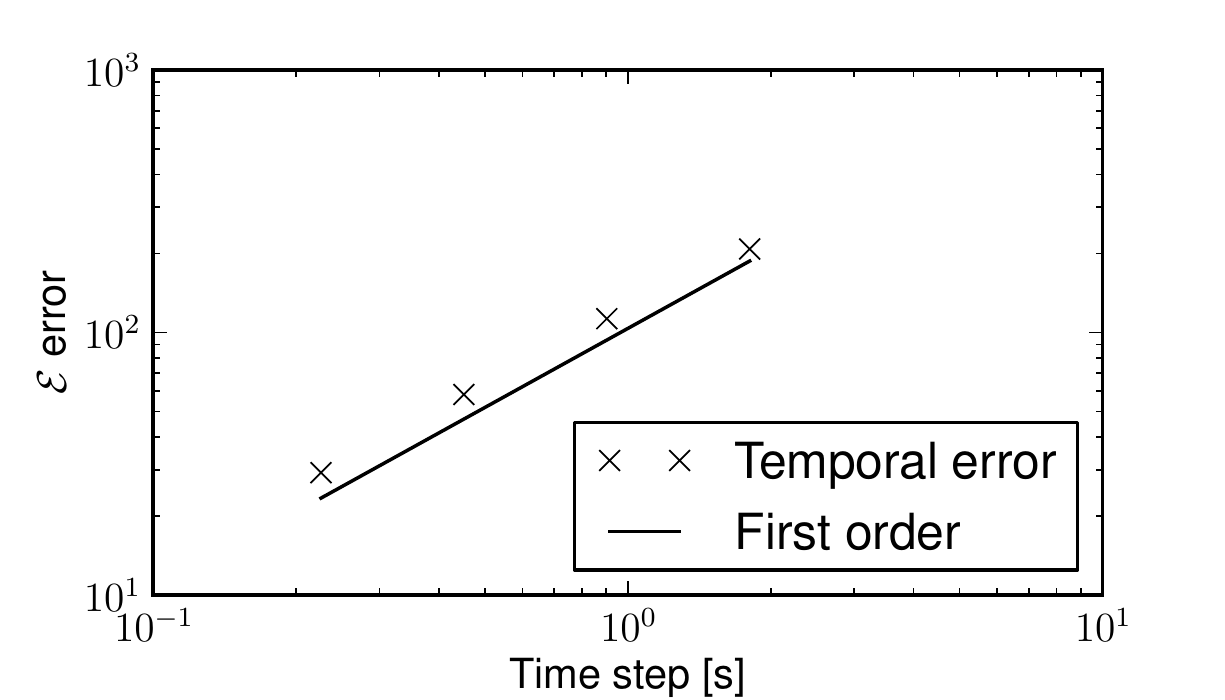}
\label{fig:results/mms_scenario/temporal_rate_of_convergence} 
}
   \caption{The expected and achieved orders of convergence for the forward model.}\label{fig:results/mms_scenario/temporal_and_spatial_rate_of_convergence} 
\end{figure}

To determine the spatial order of convergence of the model, the time step was fixed to a small value of $\Delta t
= T/16800$~s, to ensure that the numerical error of the spatial discretisation dominates the
overall discretisation error.  Then, the forward model with the added source term was run on four
uniform, increasingly fine meshes and the error in the numerical solution $(u, \eta)$ measured as:
\begin{equation*}
 \mathcal E = \left(\left|\left|u -  u_{\mathrm{exact}}\right|\right|^2_{\Omega \times (0, T)} +
 \left|\left| \eta - \eta_{\mathrm{exact}}\right|\right|^2_{\Omega \times (0, T)}\right)^{\frac{1}{2}}.
\end{equation*}
The resulting errors plotted in figure~\ref{fig:results/mms_scenario/spatial_rate_of_convergence}
show the second-order convergence that is expected from the Taylor-Hood finite element pair.

For determining the temporal order of convergence, a mesh with $2.5$~m element size in the
$x$-direction and $160$~m element size in the $y$-direction was generated (the analytical solution
does not vary in the $y$-direction and hence a relatively large mesh element size can be used).
This fine mesh resolution ensures that the numerical error of the temporal discretisation dominates the
overall discretisation error.  The forward model was then run with a set of different time steps.
The resulting errors plotted in figure~\ref{fig:results/mms_scenario/temporal_rate_of_convergence}
show the first-order convergence expected for the implicit Euler time discretisation.

\subsection{Verification of the gradient computation}
The gradient computation was verified using the Taylor remainder convergence test.
Let $\hat{J}(m) \equiv J(z(m), m)$.
The first-order Taylor expansion states that:
\begin{equation}
  \left|\hat{J}(m +  {\delta m}) - \hat{J}(m) - \frac{d\hat{J}}{dm}{\delta m}\right| = \mathcal{O}(\|{\delta m}\|^2). 
\end{equation}
Examining the convergence order of the remainder term is a strong test that the adjoint model and the gradient evaluation are implemented correctly:
for nonlinear functionals, the gradient computed using the discrete adjoint approach is correct if and only if the Taylor
remainder converges at second order.

Firstly, a simple configuration with a single turbine was set up with the turbine parameterisation and the functional of interest as described above.
The Taylor remainder convergence test in many random perturbation directions $\delta m$ yielded the expected second-order convergence (not shown).
Secondly, the Taylor remainder convergence test was applied on several of the numerical examples
presented in the following section (not shown). All yielded second-order convergence, giving confidence that the
adjoint model and the gradient computation are implemented correctly.
 
\section{Examples}\label{sec:examples}
\begin{figure}[p]
  \centering
    \subfloat[Scenario 1]{
      \centering
      \includegraphics[width=0.39\textwidth]{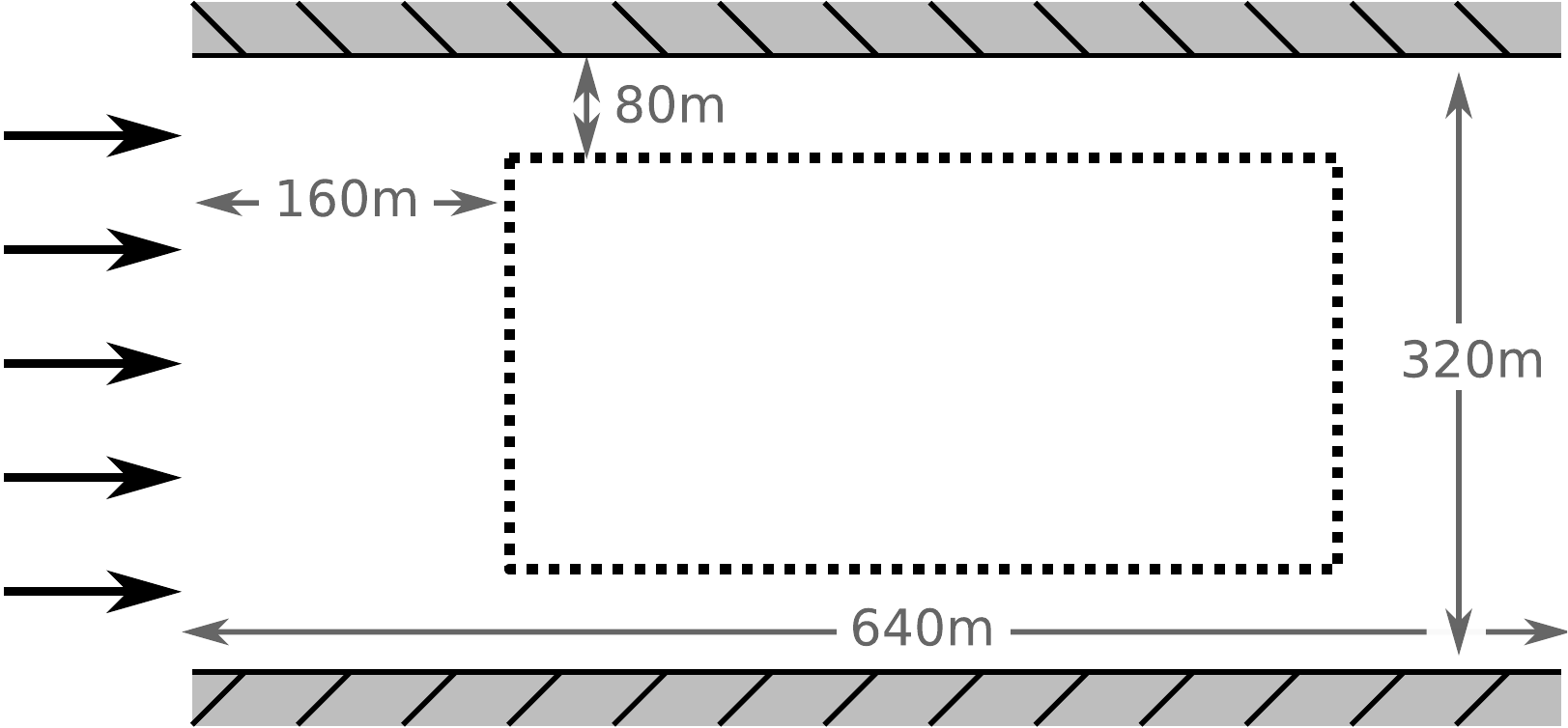}
      \label{fig:scenarios/scenario1}
    }
    \subfloat[Scenario 2]{
      \centering
      \includegraphics[width=0.39\textwidth]{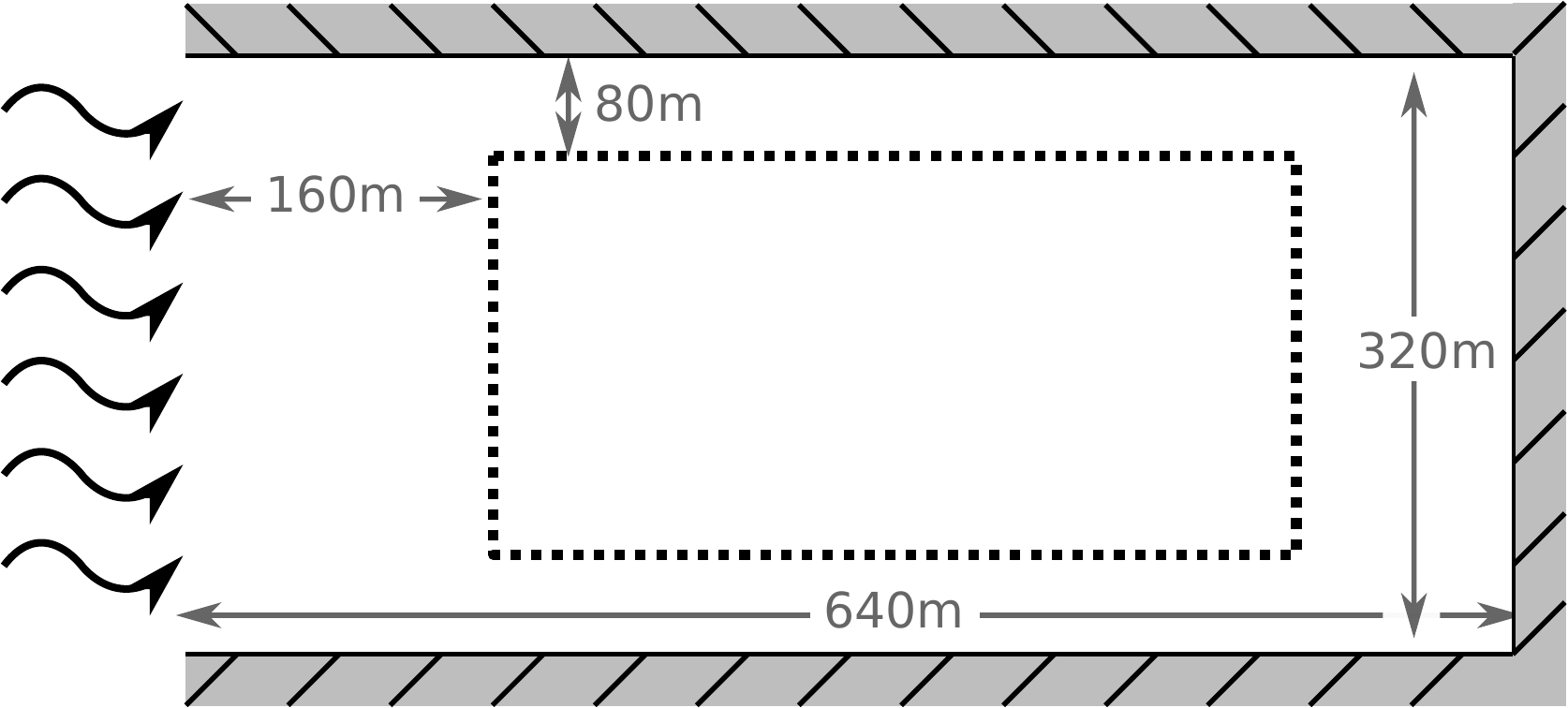}
      \label{fig:scenarios/scenario2}
    }
    \\
    \subfloat[Scenario 3]{
      \centering
      \includegraphics[width=0.3\textwidth]{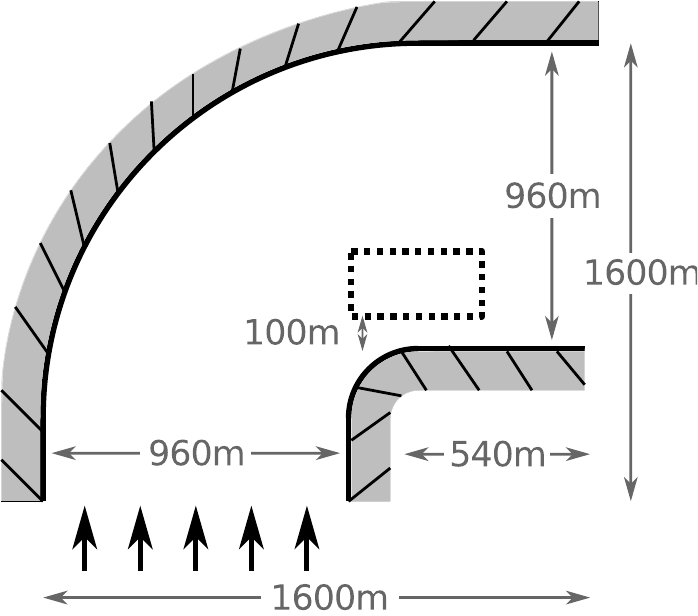}
      \label{fig:scenarios/scenario3}
    }
    \hspace{1cm}
    \subfloat[Scenario 4]{
      \centering
      \includegraphics[width=0.3\textwidth]{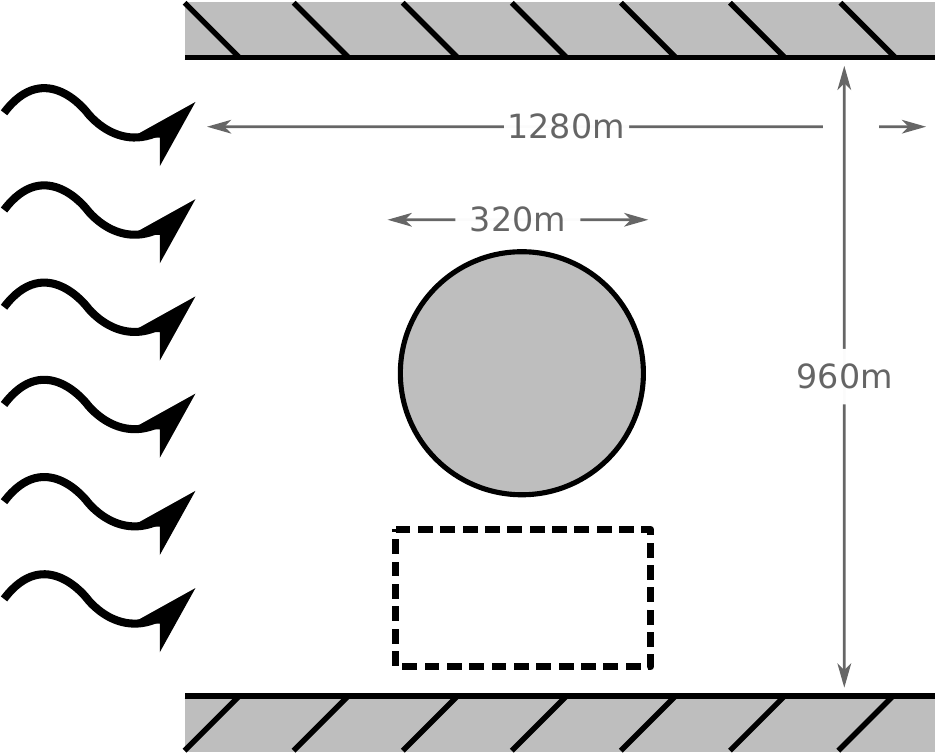}
      \label{fig:scenarios/scenario4}
    }
    \caption{The four turbine optimisation scenarios considered in the numerical examples, motivated by \cite{draper2009}. The dashed lines mark the $320\textrm{ m} \times 160\textrm{ m}$ sized sites where $32$ turbines are to be deployed. In scenarios 1 and 3 a constant inflow velocity is enforced, while the other scenarios are driven by a sinusoidal inflow.}\label{fig:draper2009}
\end{figure}

The following numerical examples solve the optimal layout problem in four idealised scenarios motivated by \cite{draper2009} (figure~\ref{fig:draper2009}). 
Additionally, in scenario 3 we optimise for the positions and the tuning of the individual turbines. 
In all examples, $32$ turbines are to be deployed in a rectangular turbine site of size $320\textrm{ m} \times 160\textrm{ m}$.
The idealised domains simplify the subsequent interpretation of the optimised configurations.
Nevertheless, the domains are chosen such that they resemble structures that can be found in practical deployment sites.
The main three objectives for the numerical examples are to investigate: 
by how much can optimisation increase the energy extraction?
Can the optimisation algorithm reliably improve the energy extraction for different scenarios?
How does the choice of the optimisation variables and the constraints impact the resulting configuration?

\begin{table}[p]
 \centering
\begin{tabular}{l | l}
 \toprule
 Parameter & Value \\
 \midrule
 Water depth & $H = 50$ m \\
 Viscosity coefficient & $\nu = 3\textrm{ m}^2\textrm{s}^{-1}$ \\
 Turbine friction coefficient & $K = 21$ \\
 Acceleration due to gravity & $g = 9.81\textrm{ m}\textrm{s}^{-2}$ \\
 Water density & $\rho = 1, 000 \textrm{ kg}\textrm{m}^{-3}$ \\
 Bottom friction coefficient & $c_b = 0.0025$ \\
 Turbine radii & $r=10$~m \\
 \bottomrule
\end{tabular}
\caption{The parameter values used in the experiments of section~\ref{sec:examples}. The non-dimensional bottom friction coefficient is a common value for coastal modelling \citep{vennell2012b}.}
 \label{tab:computational_settings}
\end{table}

The parameter choices for the experiments are listed in table~\ref{tab:computational_settings}. 
The stopping accuracy of the SLSQP optimisation algorithm was set to $10^{-6}$ unless
otherwise stated. The stopping criteria demand that the solution is feasible (i.e. that the $L^1$ norm of the constraint residuals is less than the tolerance) 
and that the solution is optimal (i.e. that the gradient in the search direction is also less than the tolerance).
This tolerance is extremely tight, as can be seen in the convergence plots, and could be
weakened for efficiency reasons.
Scenarios 1 and 3 were modelled using the stationary shallow water equations with the following boundary conditions.
On the inflow boundary a constant inflow velocity of $2\textrm{ m} \textrm{s}^{-1}$ was enforced and on the outflow boundary the free-surface displacement was set to zero.
A no-normal flow boundary condition with a free-slip condition for the tangential components was applied on the remaining boundaries. 
Scenarios 2 and 4 solved the non-stationary shallow water equations with boundary conditions explained in the associated sections. 


All examples used unstructured meshes with a uniform mesh element size of $h=20$~m outside the site area.
Inside the site area, the mesh was structured with an element size of $h=2$~m. 
The higher resolution in the turbine site ensures that each individual turbine is well resolved, independent of its location within the site; this
obviates the need for regridding when the turbines are moved.
Doubling the resolution for the problem considered in section~\ref{sec:turb_opt_a_single_turbine} changed the power extracted by less than $0.5\%$.
It is therefore assumed that the problems are sufficiently well resolved.
The resulting meshes consisted of approximately $33, 000$ triangles for scenario 1 and 2, $63, 000$ triangles for scenario 3 and $45, 000$ triangles for scenario 4.
All meshes were generated using Gmsh \citep{geuzaine2009}.

In all numerical experiments, the optimisation algorithm was initialised with the $32$ turbines deployed in a regular $8 \times 4$ grid
and with box constraints for the turbine positions to ensure that the turbines remain inside the site areas.

\subsection{A single turbine}\label{sec:turb_opt_a_single_turbine}
\begin{figure}
    \centering
    \subfloat[Turbine position]{
      \centering
      \includegraphics[width=0.39\textwidth]{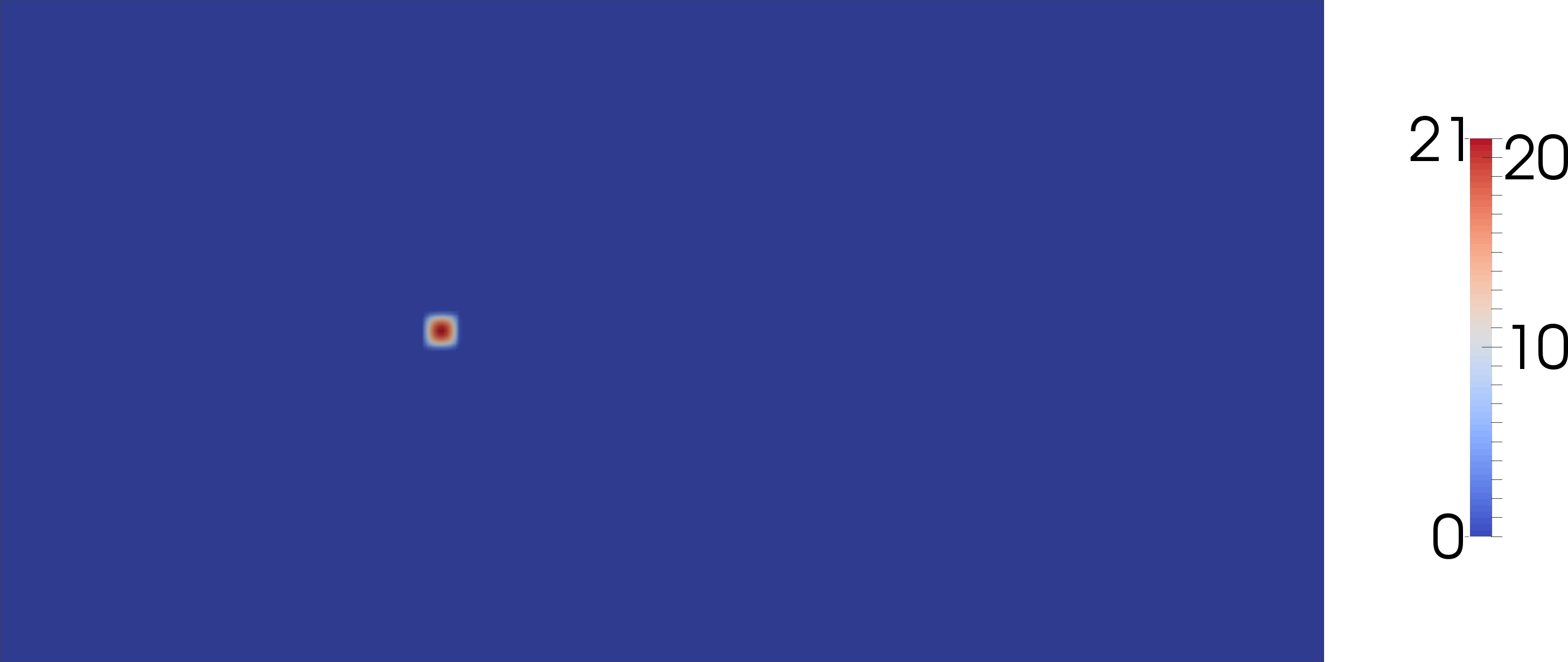}
      \label{fig:results/scenario0/paraview/turbine_friction}
    }
    \subfloat[Velocity magnitude]{
      \centering
      \includegraphics[width=0.39\textwidth]{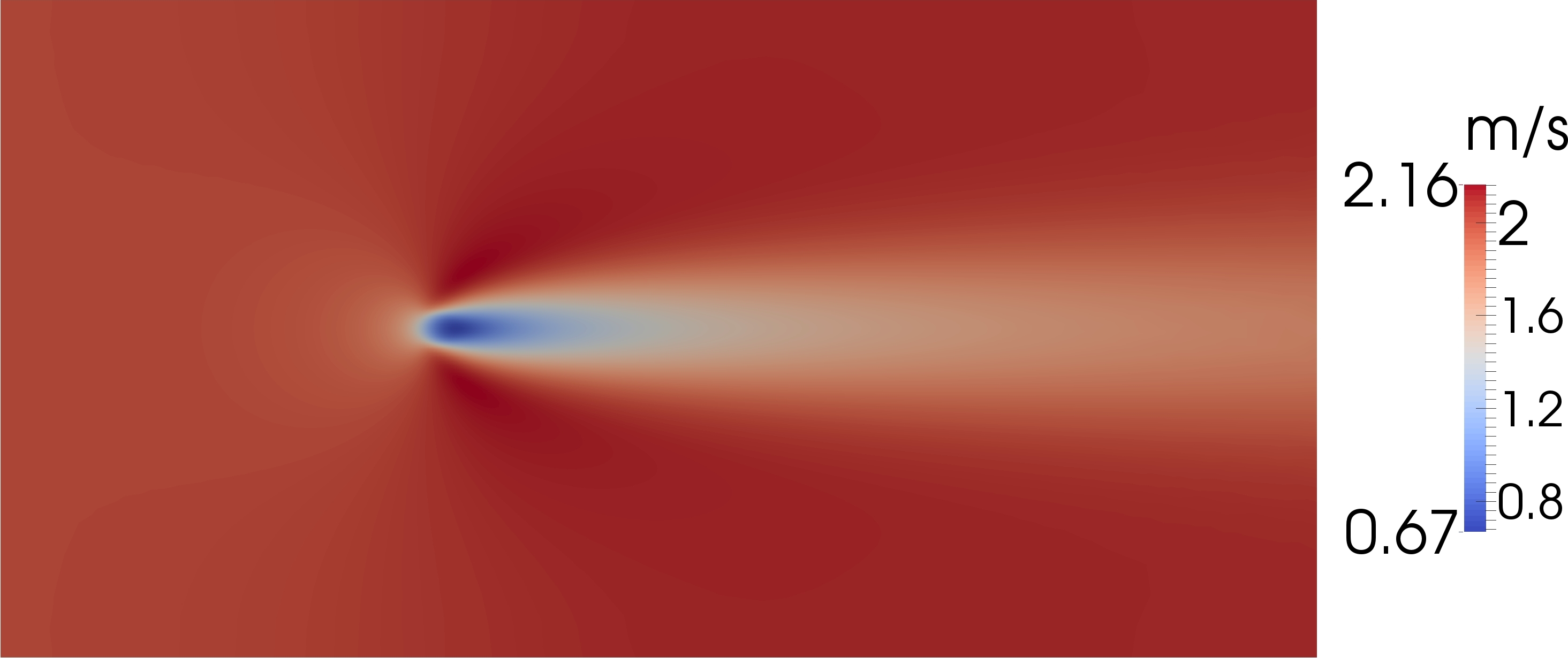}
     }
     \\
    \subfloat[Dependency of the power extraction $J$ on the friction coefficient $K$ in equation~\eqref{eq:single_turbine_coefficient}]{
       \centering
       \includegraphics[width=0.39\textwidth]{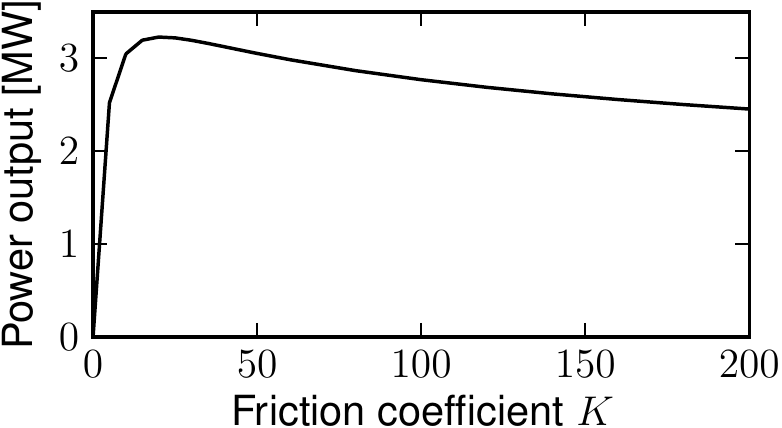}
       \label{fig:results/scenario0/turbine_friction_vs_power} 
     }
  \caption{The results of deploying a single turbine in the domain of scenario 1.} 
\end{figure}

As a preliminary test, a single turbine was deployed in the setup of scenario 1.
The turbine was placed $640/3\textrm{~m}\times320/2\textrm{~m}$ away from the bottom left corner, as shown in figure~\ref{fig:results/scenario0/paraview/turbine_friction}.

This setup was used to study the dependency of the power extraction $J$ on the friction coefficient $K$ that occurs in the turbine parameterisation~\eqref{eq:single_turbine_coefficient}. 
Figure~\ref{fig:results/scenario0/turbine_friction_vs_power} shows the power extraction for a range of $K$ values.
The graph shows a defined single peak where the power extraction is maximised, and is similar to previous studies~\citep{garrett2005, vennell2012, vennell2012b}. 
The reason for this peak is that as $K \to 0$, the turbine friction function $c_t$ approaches $0$, which in turn results in no power extraction since the power function~\eqref{eq:functional_steady} is multiplied by $c_t$;
similarly, as $K \to \infty$, the flow is deflected around the turbine and results in the observed power drop. 
The power extraction peaks for $K = 21$, which was used for all following numerical tests if not otherwise stated.
With this choice the single turbine extracted $3.2$~MW from the flow. 

\subsection{Scenario 1}
\begin{figure}[p]
  \centering
\subfloat[Initial turbine positions]{
    \centering
      \includegraphics[width=0.39\textwidth]{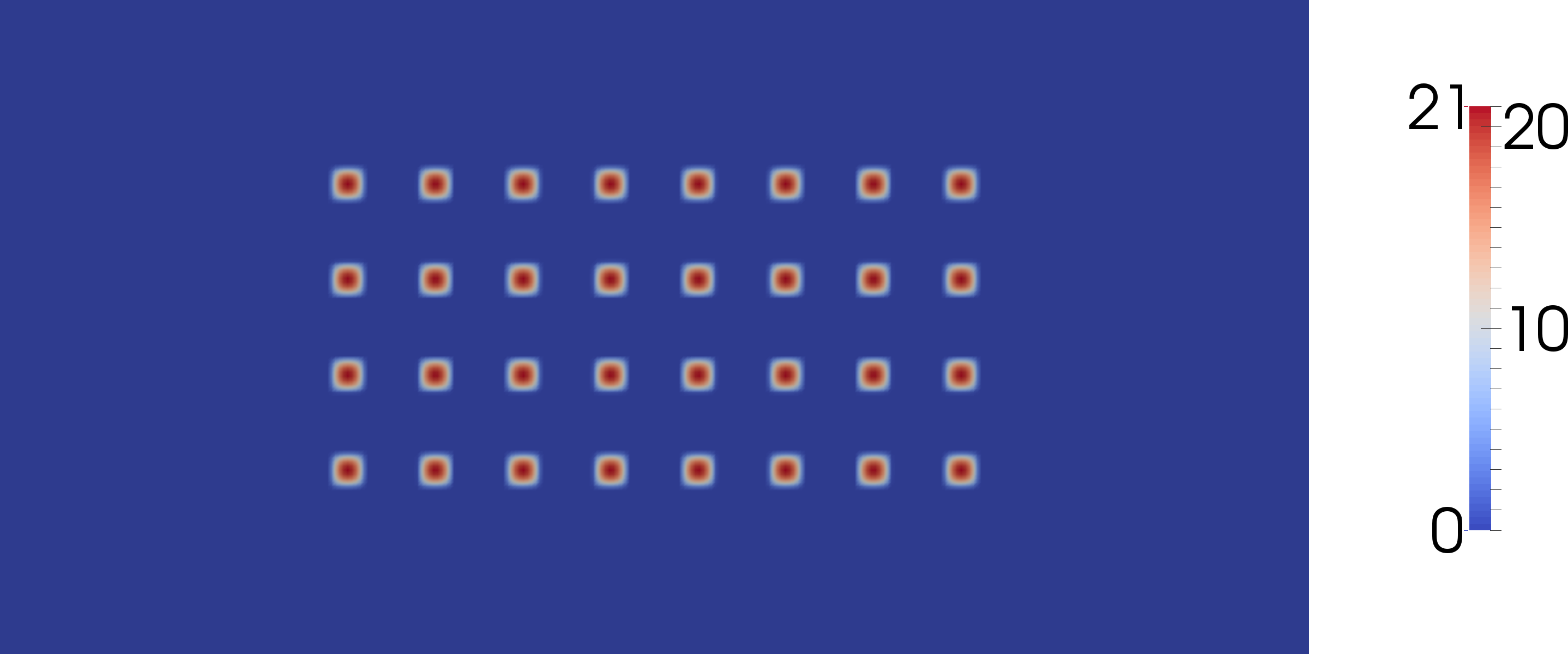}
\label{fig:results/scenario1_no_ineq/paraview/initial_turbine_friction}
}
    \subfloat[Optimised turbine positions]{
    \centering
      \includegraphics[width=0.39\textwidth]{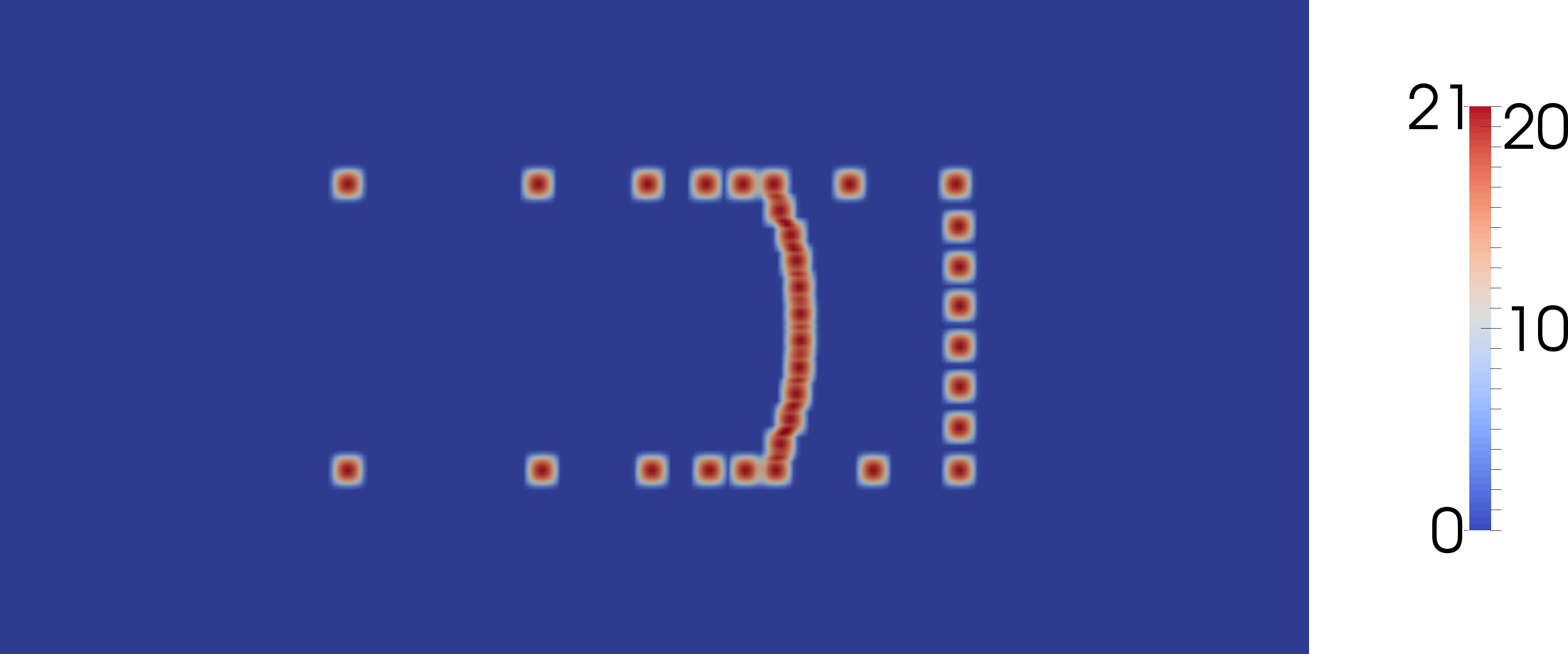}
\label{fig:results/scenario1_no_ineq/paraview/turbine_friction}
}
\\
    \centering
    \sbox{\tempbox}{
     \hbox{\includegraphics[width=0.39\textwidth]{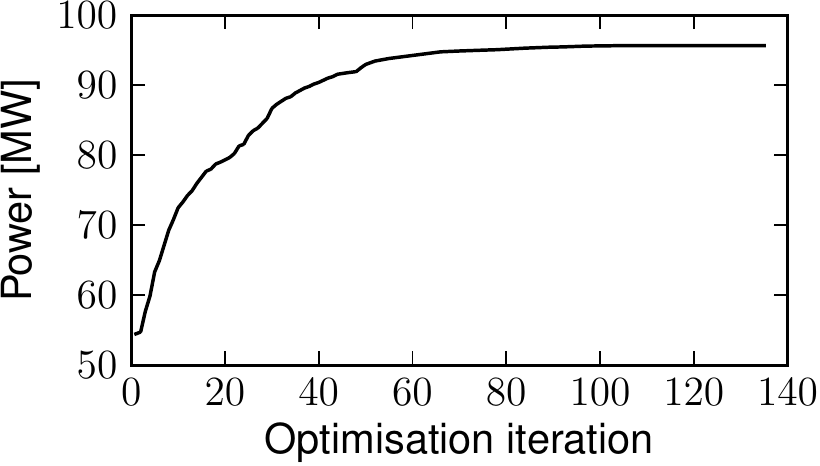}
     }
}
\subfloat[Optimisation convergence]{\usebox{\tempbox}
\label{fig:results/scenario1_no_ineq/iter_plot} 
}
    \subfloat[Velocity magnitude]{
     \vbox to \ht\tempbox{
      \vfill
      \hbox{\includegraphics[width=0.39\textwidth]{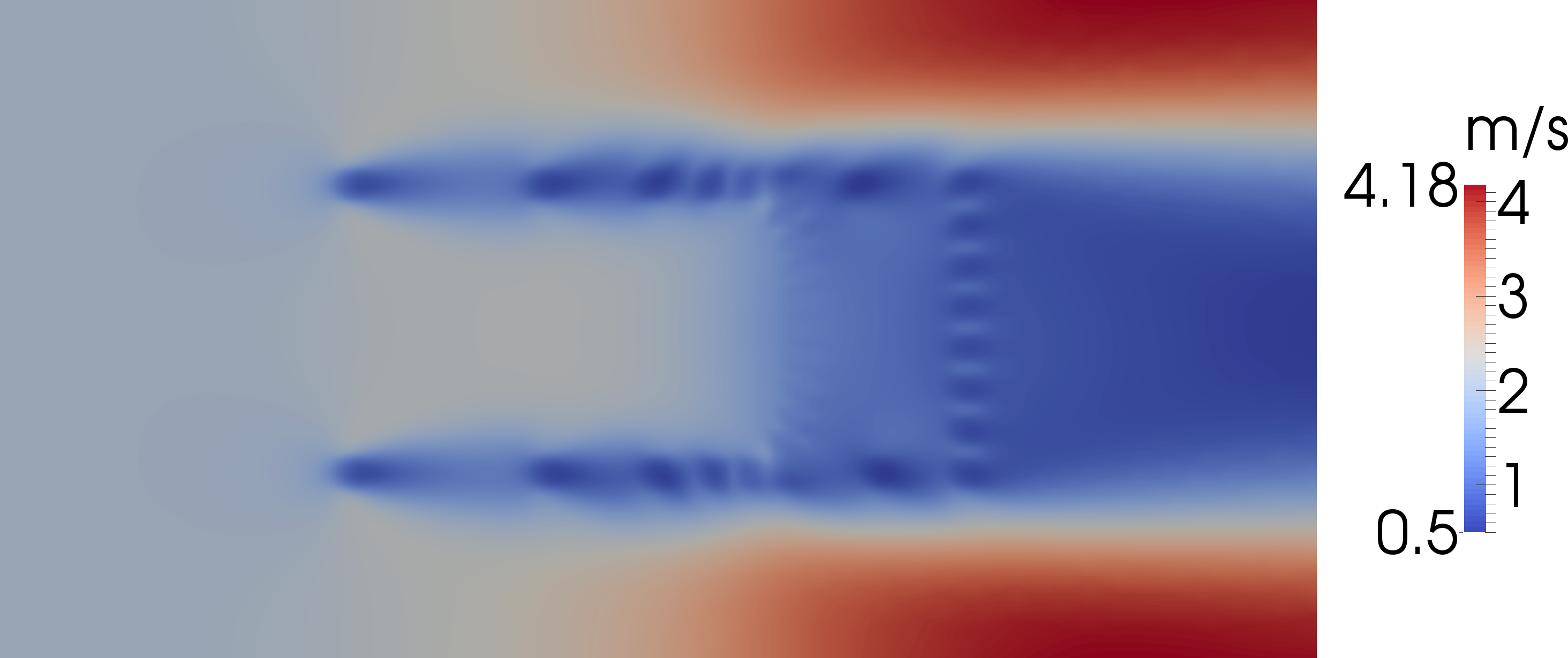}
      }
      \vfill
\label{fig:results/scenario1_no_ineq/paraview/velocity}
     }
}
\caption{Results of scenario 1 without minimum distance constraints, i.e. turbines may overlap.}\label{fig:results/scenario1_no_ineq} 
\end{figure}

Firstly, the layout problem for scenario 1 (figure~\ref{fig:scenarios/scenario1}) was solved without enforcing a minimum distance between turbines, 
i.e. the turbines can be placed arbitrarily inside the site area and may even overlap.
With that setup the optimisation algorithm terminated after $135$ iterations ($134$ gradient evaluations,
$231$ functional evaluations).
The results are shown in figure~\ref{fig:results/scenario1_no_ineq}.
Compared to the initial regular layout (figure~\ref{fig:results/scenario1_no_ineq/paraview/initial_turbine_friction}) 
the optimisation algorithm was able to increase the farm power extraction by $76\%$ from $54.5$~MW to $95.7$~MW (figure~\ref{fig:results/scenario1_no_ineq/iter_plot}).
The optimised layout (figure~\ref{fig:results/scenario1_no_ineq/paraview/turbine_friction}) has the turbines aligned in the shape of two $\sqsupset$s with the open end facing the inflow. 
An intuitive interpretation of this layout is that the water is funneled by the two sides of the $\sqsupset$s and then forced through the dense turbine `wall' at their closed ends. 
This interpretation is confirmed by an increasing free-surface displacement (not shown) and velocity difference along the sides of the $\sqsupset$s and the large jump along their closed end (figure~\ref{fig:results/scenario1_no_ineq/paraview/velocity}). 

\begin{figure}[p]
  \centering
      \subfloat[Initial turbine positions]{
    \centering
      \includegraphics[width=0.39\textwidth]{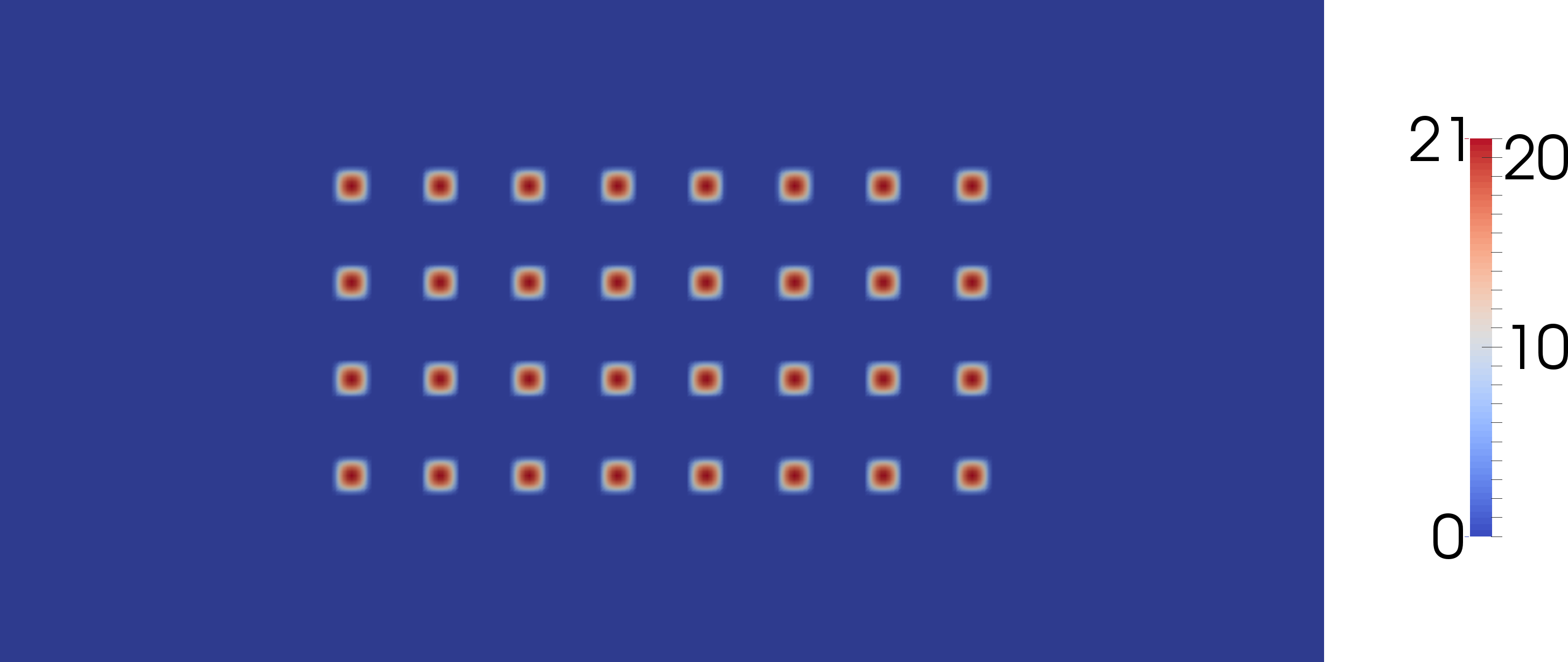}
     }
    \subfloat[Optimised turbine positions]{
    \centering
      \includegraphics[width=0.39\textwidth]{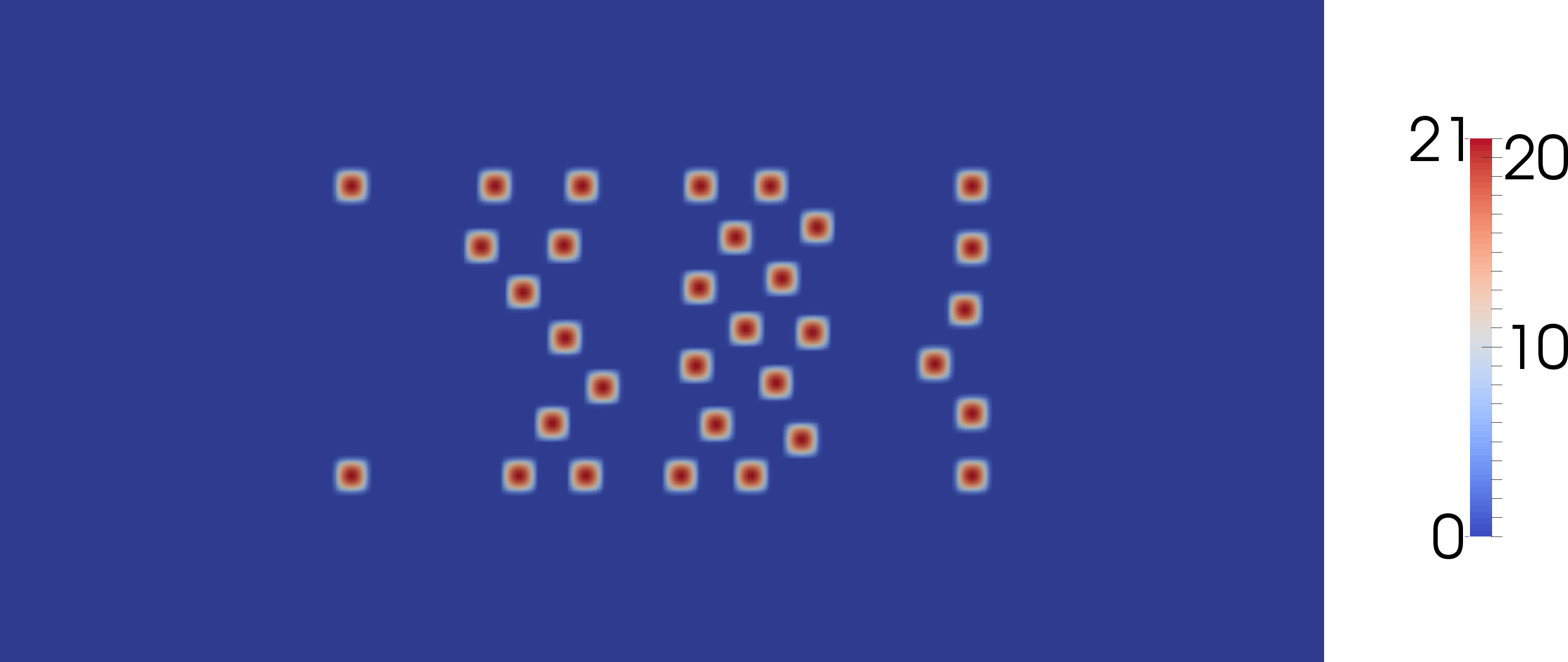}
\label{fig:results/scenario1/paraview/turbine_friction}
}
     \\
     \centering
     \sbox{\tempbox}{
      \hbox{\includegraphics[width=0.39\textwidth]{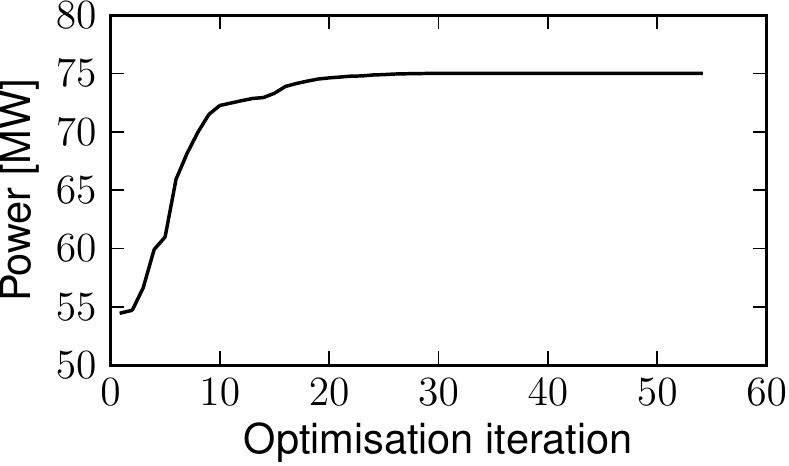}}
    }
    \subfloat[Optimisation convergence]{\usebox{\tempbox}
   \label{fig:results/scenario1/iter_plot}}
\subfloat[Velocity magnitude]{
     \vbox to \ht\tempbox{
      \vfill
      \hbox{\includegraphics[width=0.39\textwidth]{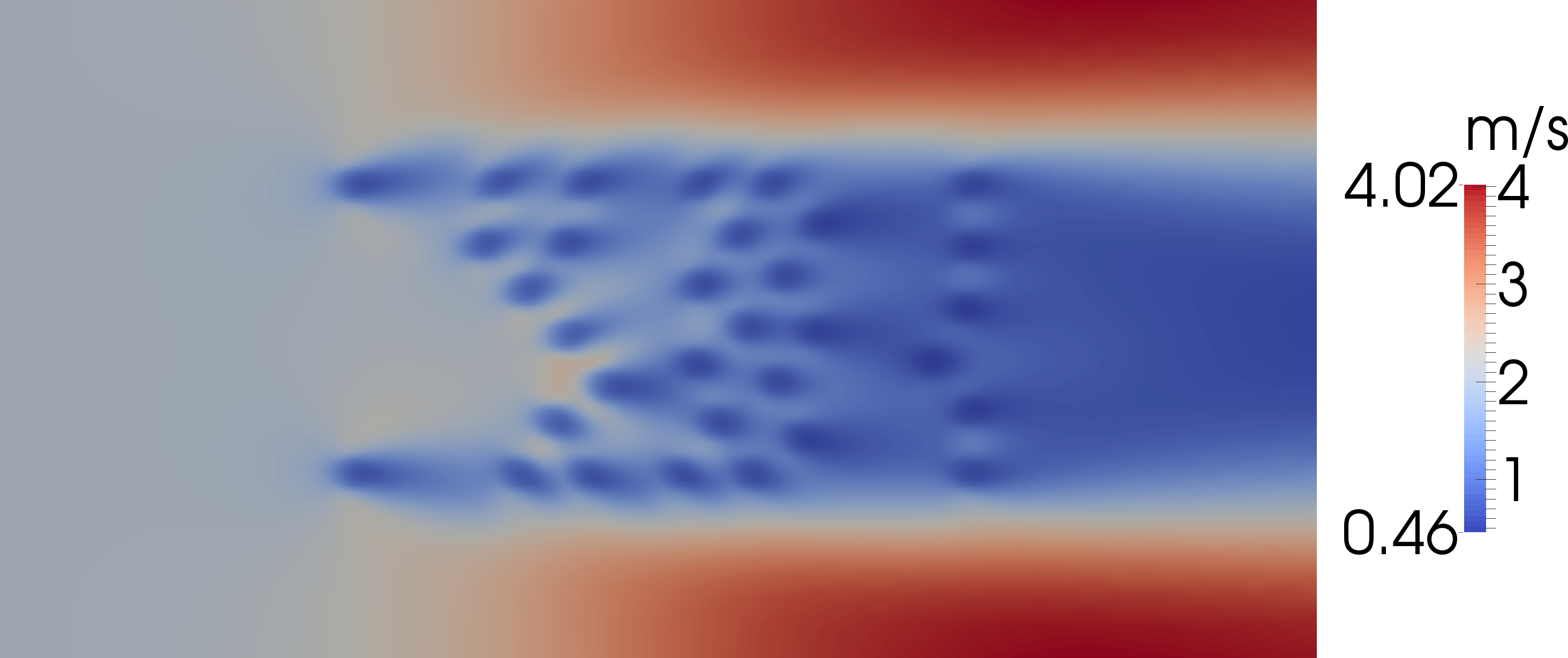}
      }
      \vfill
\label{fig:results/scenario1/paraview/velocity}
     }
}
  \caption{Results of scenario 1 with minimum distance constraints.}\label{fig:results/scenario1} 
\end{figure}

An additional experiment was performed where inequality constraints were included to enforce a minimum
distance of $30$~m ($3$ turbine radii) between each turbine.  With this setup, the optimisation algorithm
terminated after $54$ iterations ($53$ gradient evaluations, $112$ functional evaluations) and the farm power
extraction increased by $38\%$ from $54.5$ MW~to $75.0$~MW.  The reduced optimised power extraction
compared to the previous setup is expected since the inequality constraints add further restrictions to the
feasible turbine positioning.  In particular, the previous optimised turbine layout is not a feasible solution
for this setup.

The optimised alignment differs significantly from the previous one (figure \ref{fig:results/scenario1/paraview/turbine_friction}).
The two main characteristic structures are a $>$ shaped alignment close to the inflow boundary and a wall of turbines near the outflow boundary. 
Furthermore, the turbines are staggered to avoid placing one turbine in the direct wake of another turbine.
Finally, the free-surface displacement (not shown) and the velocity magnitude decrease more gradually towards the outflow compared to the previous setup (figure~\ref{fig:results/scenario1/paraview/velocity}). 

\subsection{Scenario 2}\label{sec:turbine_opt_scenario2}
The layout problem for scenario 2 (figure~\ref{fig:scenarios/scenario2}) was solved using the non-stationary shallow water equations. 

\begin{figure}
  \centering
      \subfloat[Initial turbine positions]{
    \centering
      \includegraphics[width=0.39\textwidth]{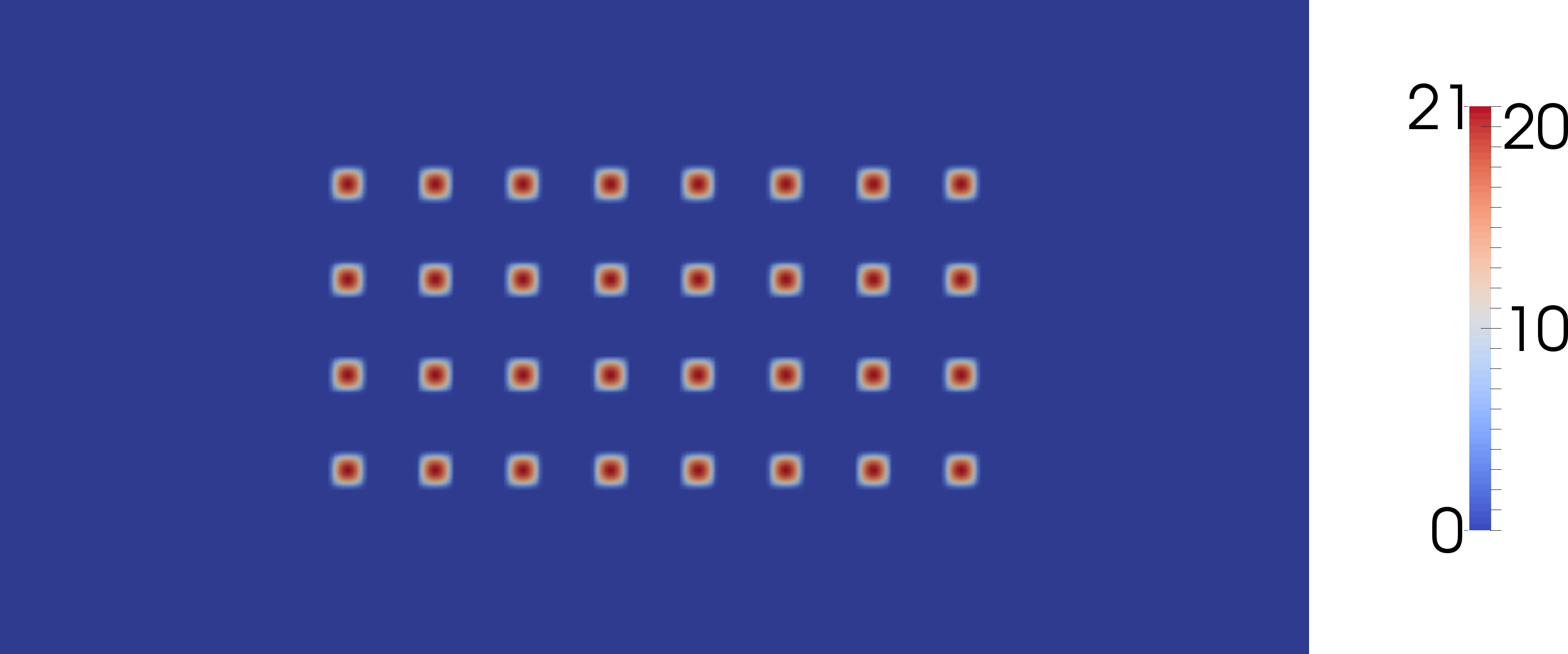}
     }
    \subfloat[Optimised turbine positions]{
    \centering
      \includegraphics[width=0.39\textwidth]{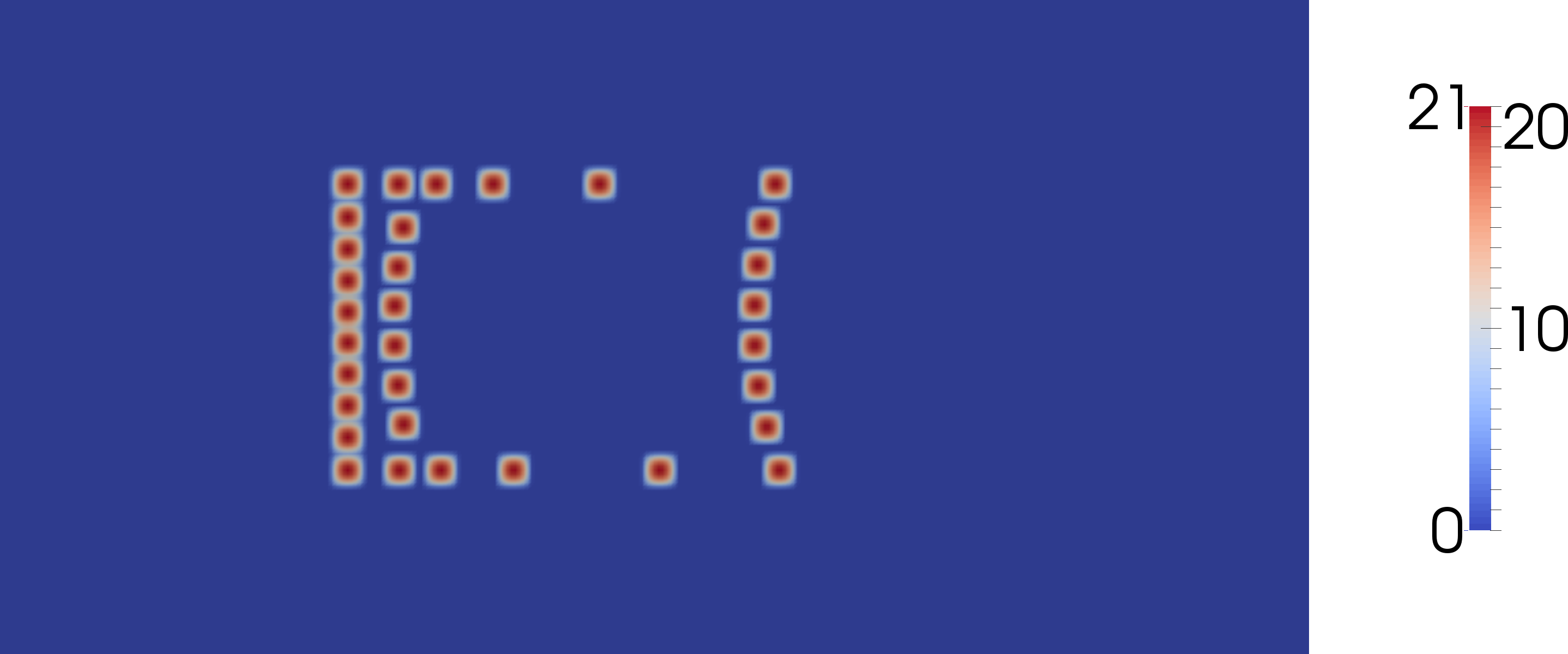}
\label{fig:results/scenario2_noineq/paraview/turbine_friction}
}
      \\
\subfloat[Velocity magnitude during the\newline ebb tide]{
    \centering
      \includegraphics[width=0.39\textwidth]{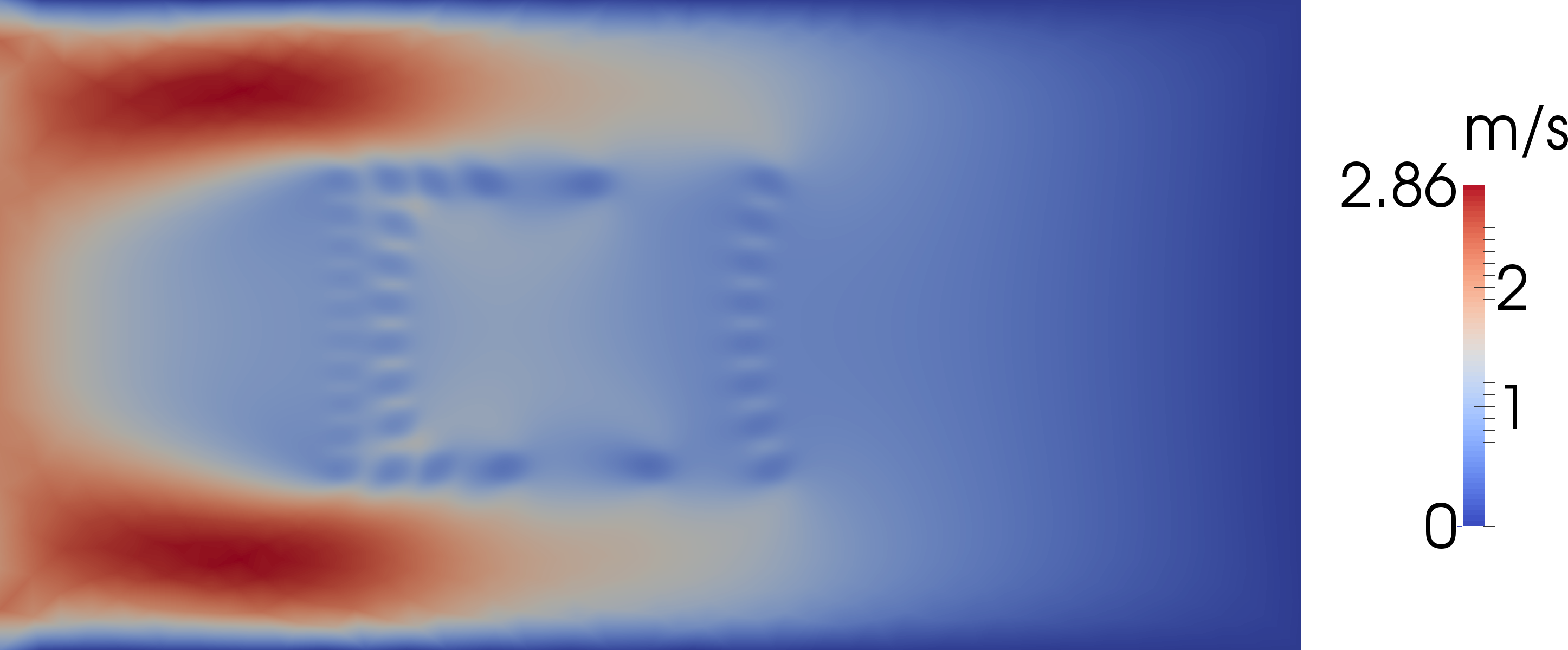}
      \label{fig:results/scenario2_noineq/paraview/velocity_ebb}
     }
\subfloat[Velocity magnitude during the\newline flood tide]{
    \centering
      \includegraphics[width=0.39\textwidth]{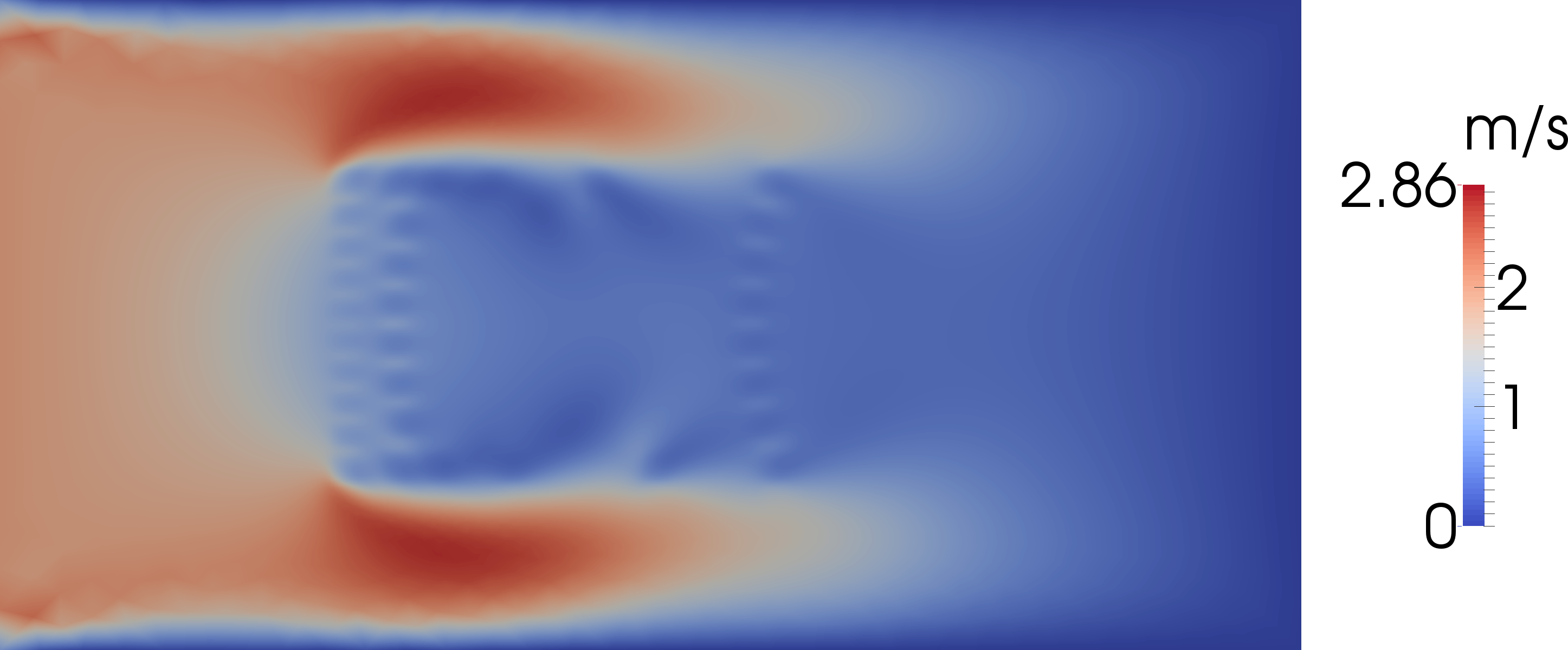}
\label{fig:results/scenario2_noineq/paraview/velocity_flood}
}
     \\
  \subfloat[Optimisation convergence]{
    \centering
    \includegraphics[width=0.39\textwidth]{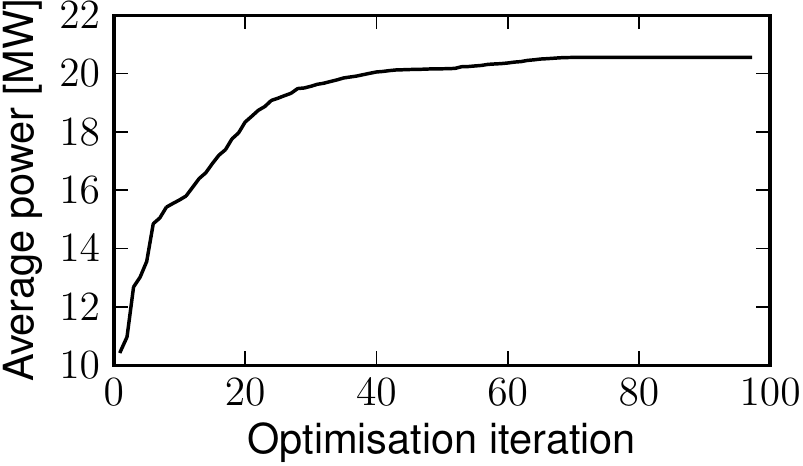}
      \label{fig:results/scenario2_noineq/iter_plot} 
      }
  \subfloat[Power extraction over time of the optimised configuration]{
    \centering
      \includegraphics[width=0.39\textwidth]{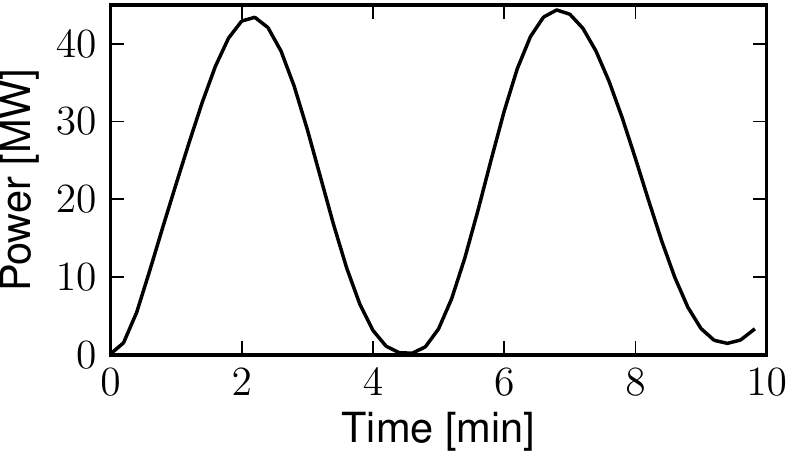}
      \label{fig:results/scenario2_noineq/power_timeline.pdf} 
      }
  \caption{Results of the non-stationary scenario 2 without minimum distance constraints.}\label{fig:results/scenario2_noineq} 
\end{figure}

The boundary conditions were as follows.
On the top, bottom and right boundaries a no-slip boundary condition was applied.
On the left boundary a Dirichlet boundary condition enforced a sinusoidal in-/outflow velocity:
\begin{equation*}
 u(x, y, t)  =
\begin{pmatrix}
-2\sin\left(2\pi t/P\right) \\
0
\end{pmatrix},
\end{equation*}
where $P$ is the tidal period time. 
Due to the small basin size, a realistically long tidal period would lead to an excessively large tidal range. 
Therefore, the tidal period was defined to be $P \equiv 10$ minutes, which resulted in a tidal range of $\pm 12$~m.
The simulation time was set to one full tidal period with a time step of $\Delta t = 12$~s.
No spin up phase was applied, as its effect is assumed to be small due to the relatively short extent of the domain.

The optimisation was performed without enforcing a minimum distance between the turbines.  After $97$
optimisation iterations ($96$ gradient evaluations, $217$ functional evaluations) the relative functional improvement in
each iteration dropped below the tolerance and the optimisation was terminated. The results are shown in
figure~\ref{fig:results/scenario2_noineq}.  The average power extracted during one tidal cycle increased by
$96\%$ from $10.5$ MW to $20.6$ MW (figure~\ref{fig:results/scenario2_noineq/iter_plot}).  The optimised
layout (figure~\ref{fig:results/scenario2_noineq/paraview/turbine_friction}) resembles the result of scenario
1 (figure~\ref{fig:results/scenario1_no_ineq/paraview/turbine_friction}), with the difference that the
opening of the $\sqsupset$ structure faces the closed basin side. 

\subsection{Scenario 3}\label{sec:turbine_opt_scenario3}
\begin{figure}[p]
  \centering
\subfloat[Initial turbine positions]{
    \centering
      \includegraphics[width=0.39\textwidth]{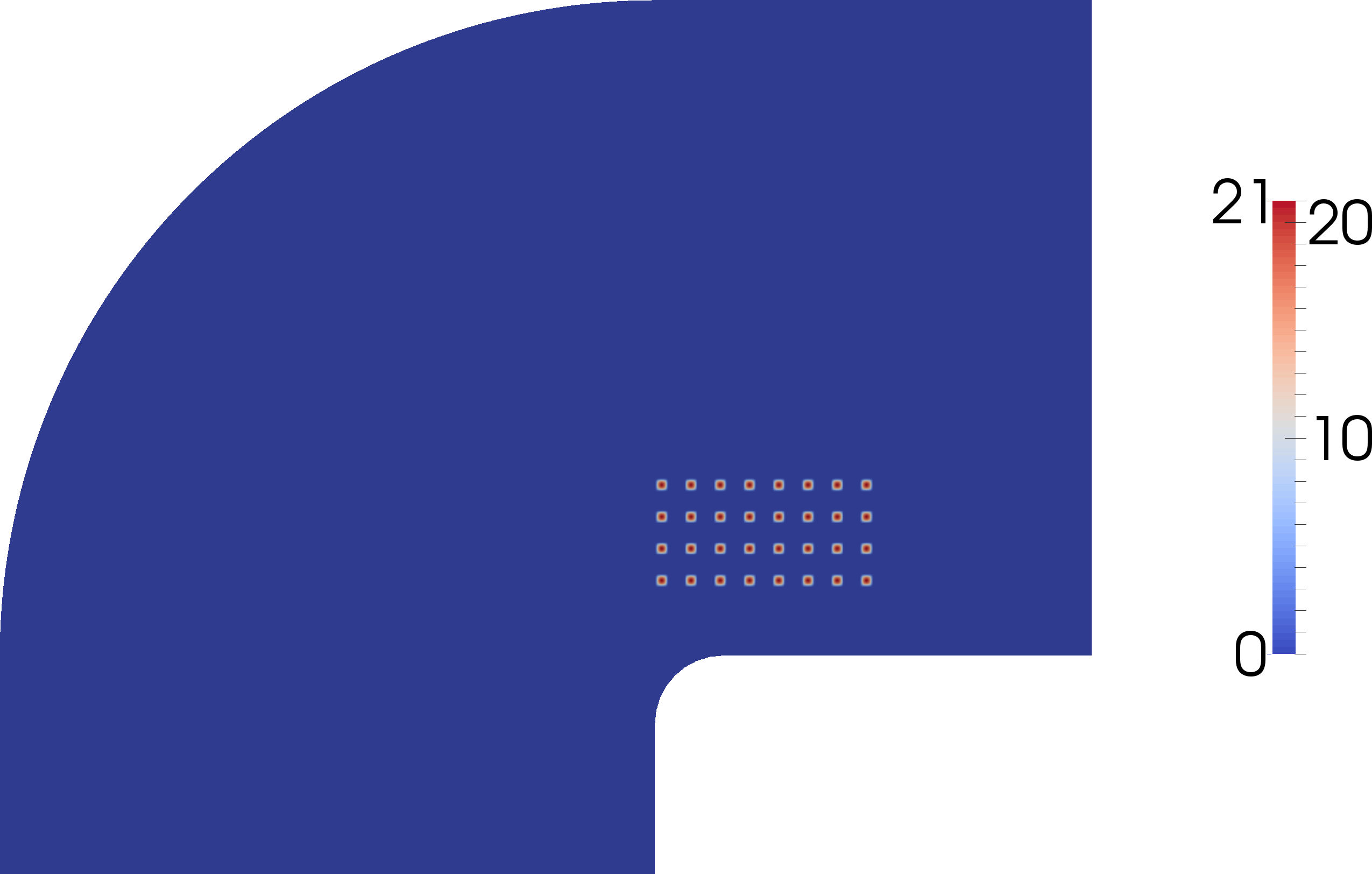}
\label{fig:results/scenario3/paraview/initial_turbine_friction}
}
\subfloat[Optimised turbine positions]{
    \centering
      \includegraphics[width=0.39\textwidth]{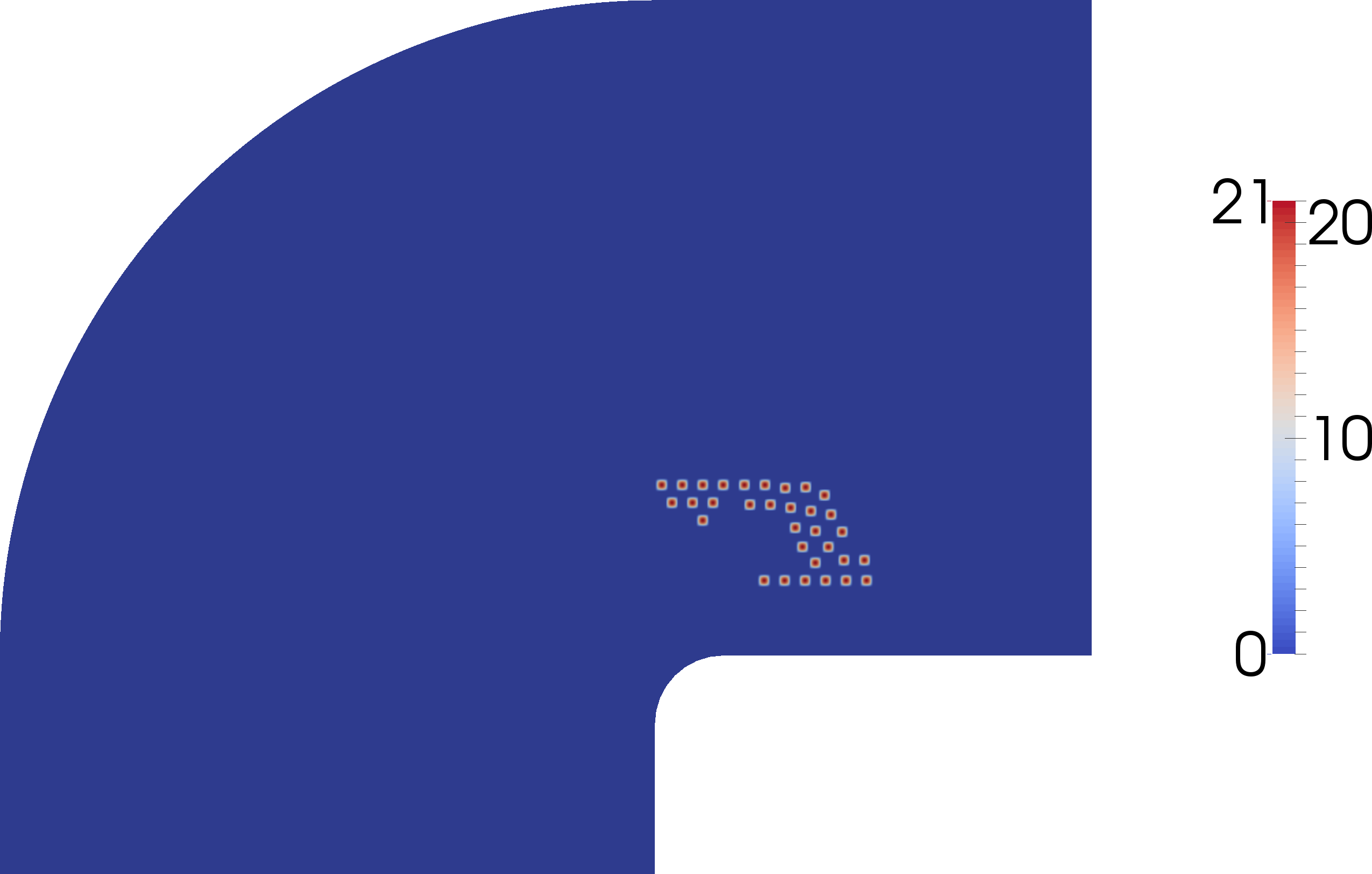}
\label{fig:results/scenario3/paraview/turbine_friction}
}
\\
\subfloat[Optimisation convergence]{
    \centering
      \includegraphics[width=0.39\textwidth]{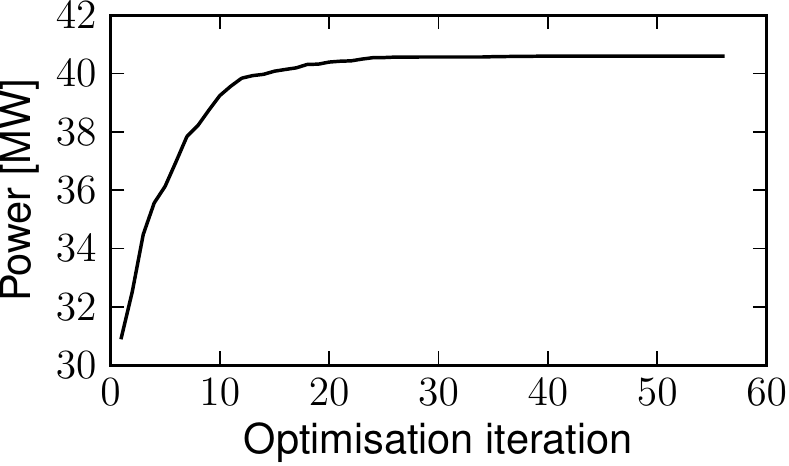}
\label{fig:results/scenario3/iter_plot} 
}
  \subfloat[Velocity magnitude]{
    \centering
      \includegraphics[width=0.39\textwidth]{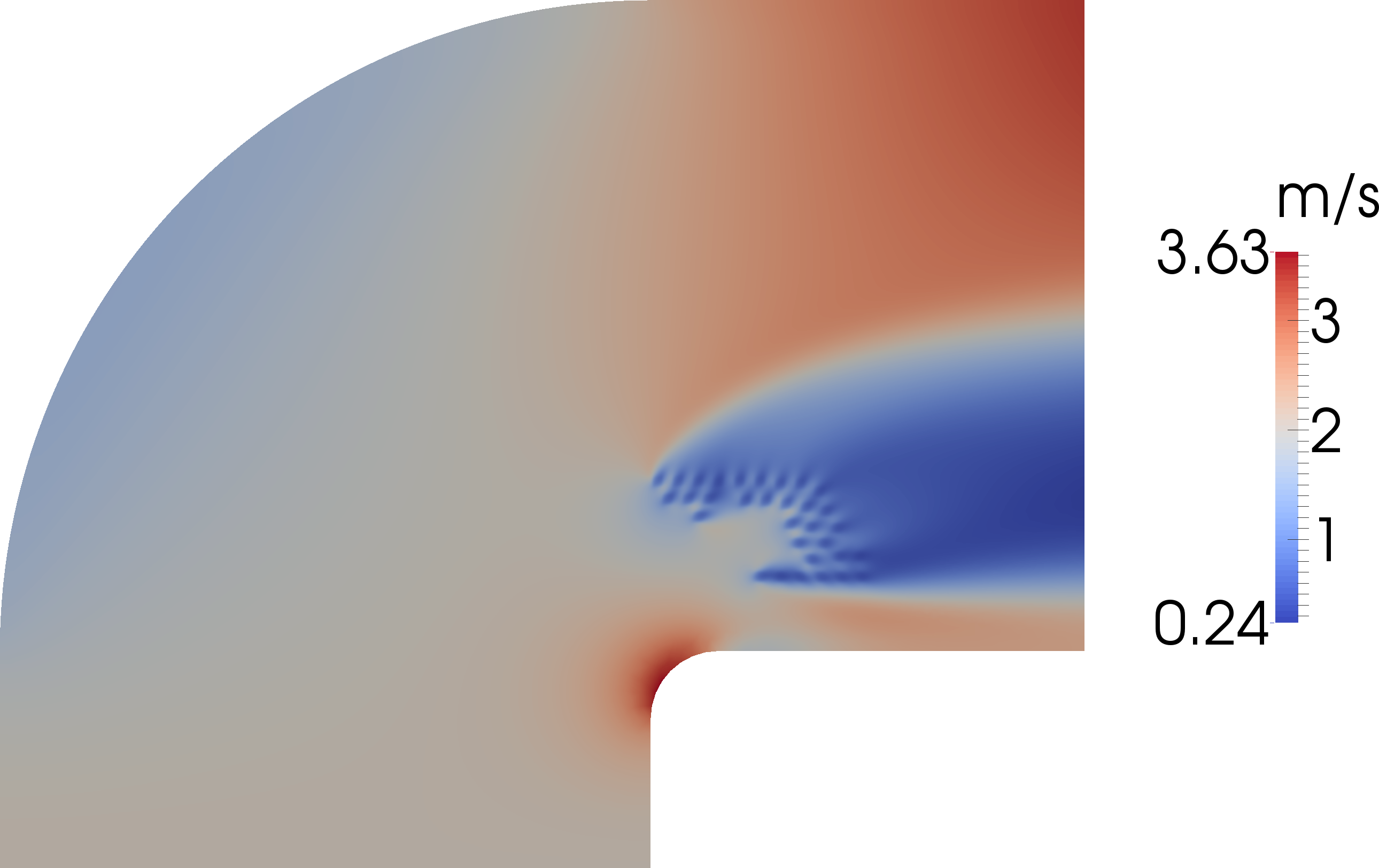}
\label{fig:results/scenario3/paraview/velocity}
}
  \caption{Results of scenario 3 with minimum distance constraints.}\label{fig:results/scenario3} 
\end{figure}

The domain of the third scenario is shown in figure~\ref{fig:scenarios/scenario3}. 
First, only the layout problem is solved.
For this test, inequality constraints were applied to enforce a minimum distance of $30$~m between each turbine.  The optimisation
algorithm terminated after $56$ iterations ($55$ gradient evaluations, $73$ functional evaluations).  The
optimised farm layout extracts $40.6$~MW, which corresponds to an increase of $31\%$ compared to the initial
layout ($30.9$~MW) (figure~\ref{fig:results/scenario3/iter_plot}).  The optimised layout features a distinct
$\Diamond-$shaped alignment with an opening on the inflow facing side
(figure~\ref{fig:results/scenario3/paraview/turbine_friction}).  Figure
\ref{fig:results/scenario3/paraview/velocity} shows the velocity magnitude and suggests that this hole acts to
trap and push the flow through the downstream turbines similar to the previous examples. 

\begin{figure}[p]
  \centering
\subfloat[Initial turbine positions and friction coefficients]{
    \centering
      \includegraphics[width=0.39\textwidth]{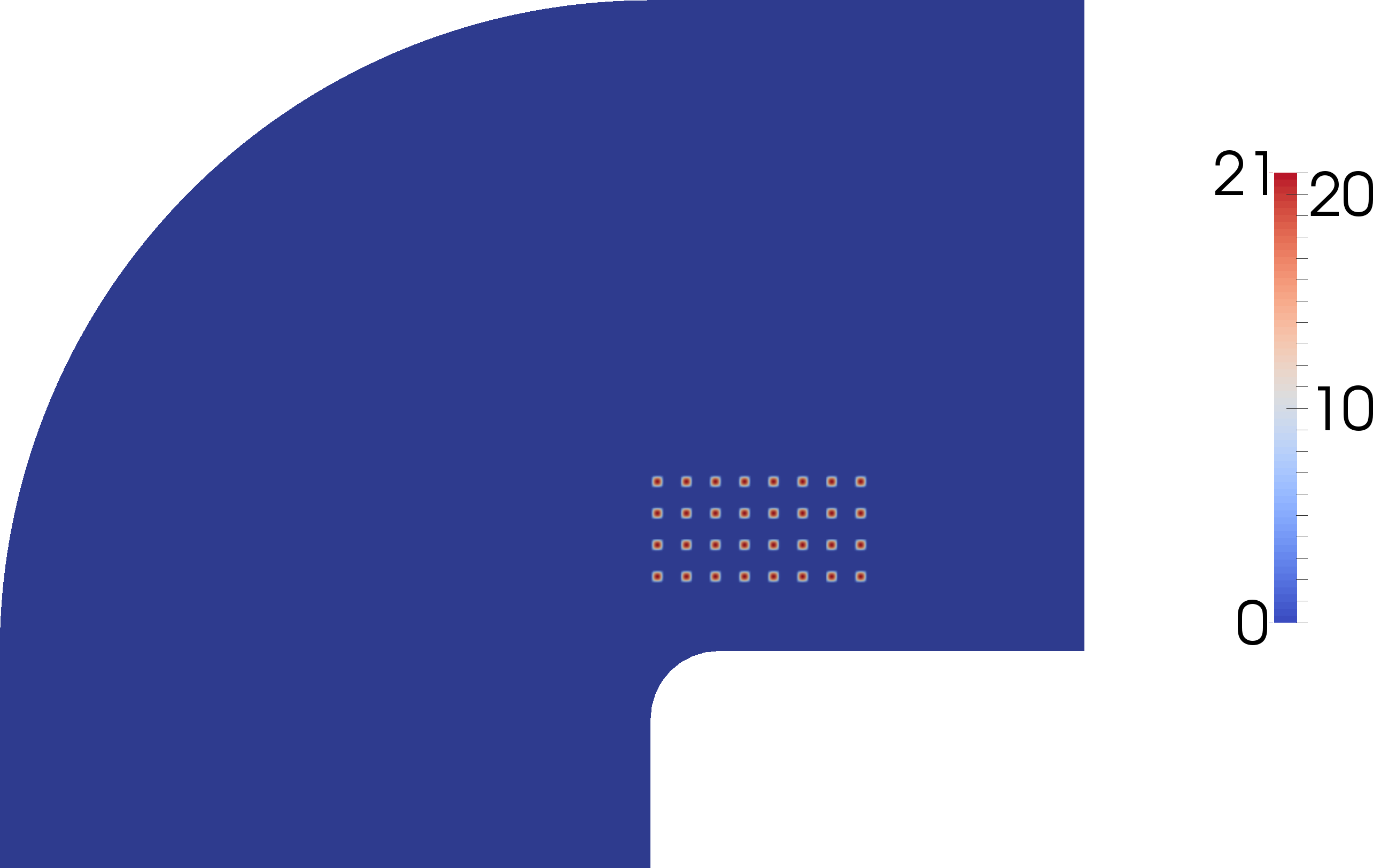}
     }
     \subfloat[Optimised turbine positions and friction coefficients]{
    \centering
      \includegraphics[width=0.39\textwidth]{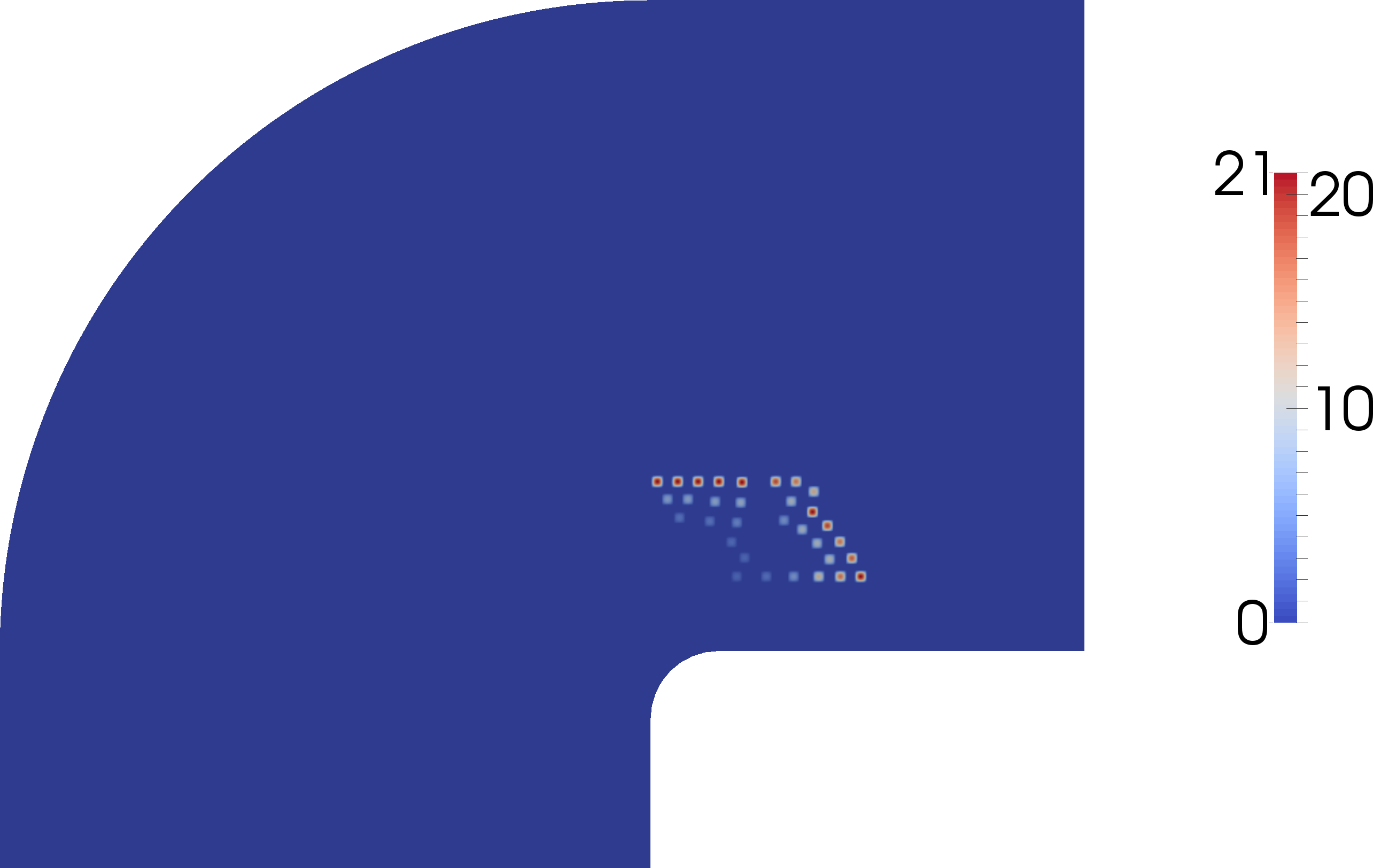}
\label{fig:results/scenario3_full_ctrl/paraview/turbine_friction}
}
     \\
  \subfloat[Optimisation convergence]{
    \centering
      \includegraphics[width=0.39\textwidth]{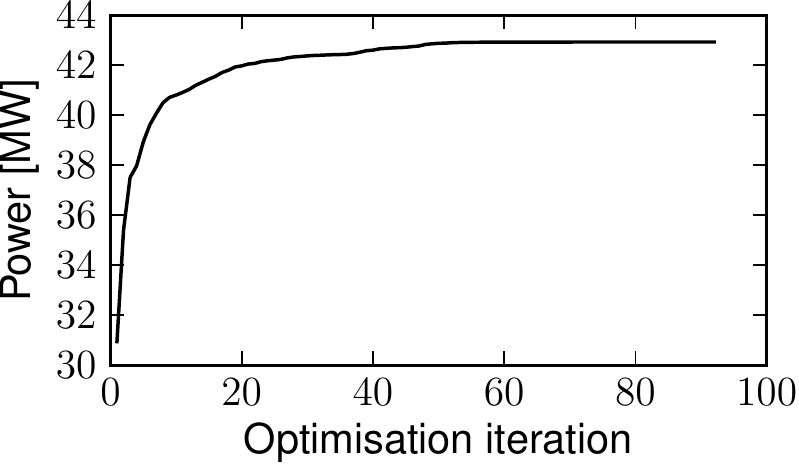}
\label{fig:results/scenario3_full_ctrl/iter_plot}
}
\subfloat[Velocity magnitude]{
    \centering
      \includegraphics[width=0.39\textwidth]{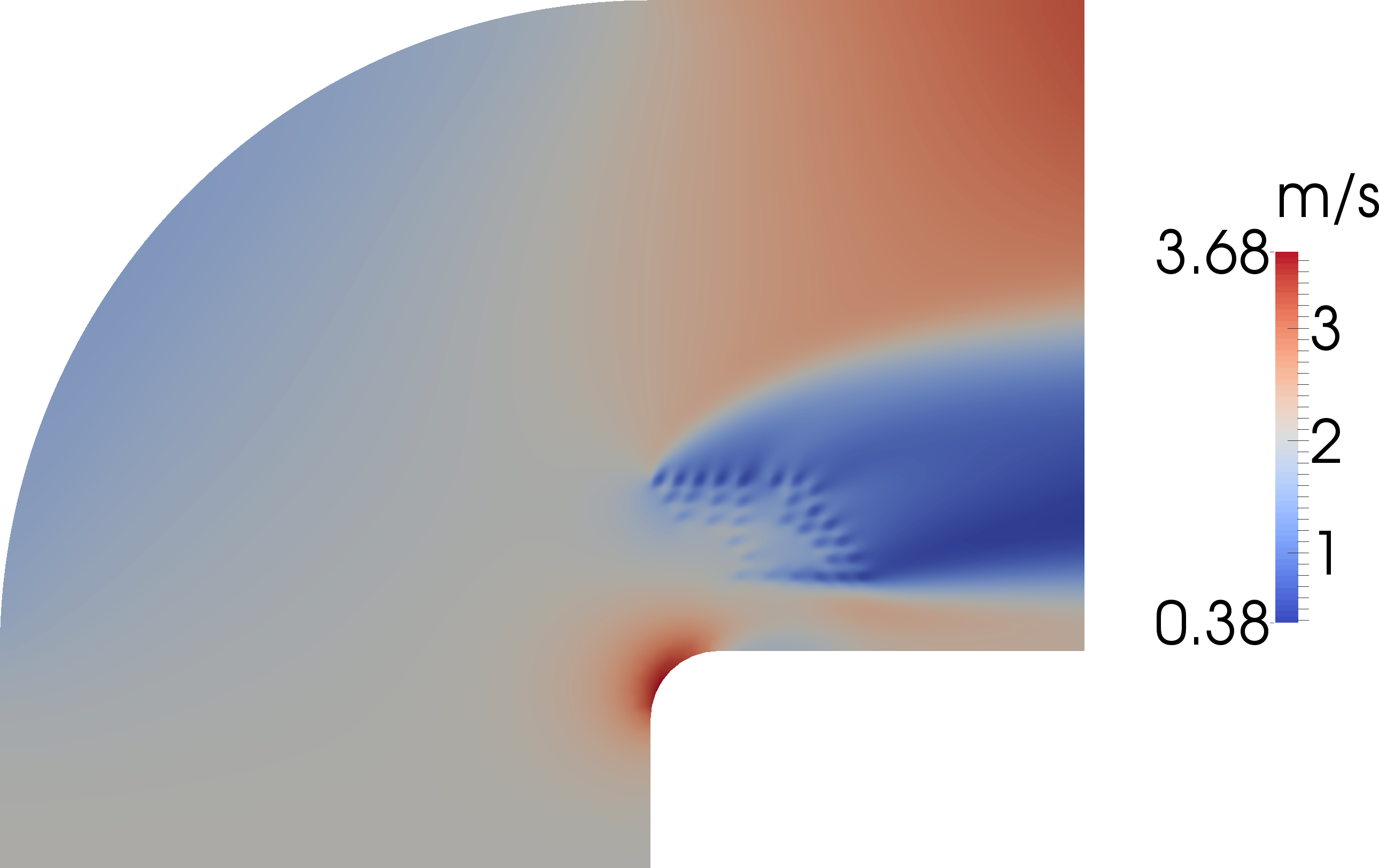}
\label{fig:results/scenario3_full_ctrl/paraview/velocity}
}
  \caption{Results of scenario 3 with minimum distance constraints and optimising for both the turbine friction coefficients and the turbine positions.}\label{fig:results/scenario3_full_ctrl} 
\end{figure}

Second, the optimisation parameters were extended to include the friction coefficients $K$ of each turbine in the parameterisation~\eqref{eq:single_turbine_coefficient}.
\citet{vennell2010, vennell2012} showed the necessity of varying the friction coefficients in order to achieve optimal farm performance. 
Each $K$ coefficient was constrained such that $0\le K \le 21$.  With this setup, the
optimisation algorithm terminated after $92$ iterations ($88$ gradient evaluations, $148$ functional
evaluations). The results are presented in figure~\ref{fig:results/scenario3_full_ctrl}. Compared to the
initial configuration, the farm power extraction increased by $39\%$ from $30.9$~MW to $42.9$~MW -- the additional
freedom of varying the $K$ coefficients resulted in a higher optimised power extraction than the previous setup.

The optimised turbine configuration is similar in shape to the previous solution but with a less distinct hole on the inflow facing side.
Most notably, the friction coefficients of most turbines are significantly reduced.
Only the turbines on the downflow edges of the $\Diamond$ take the maximum value, but nevertheless
this configuration extracts $6\%$ more energy than the optimal solution of the previous setup where only the positions were optimised.
This example shows that optimising the friction coefficient (which can be viewed as reducing the size of the turbine
or controlling the blade pitch) can lead to a significant increase in the power extraction of the farm. 

\subsection{Scenario 4}
\begin{figure}[tb]
  \centering
  \subfloat[Initial turbine positions]{
    \centering
      \includegraphics[width=0.39\textwidth]{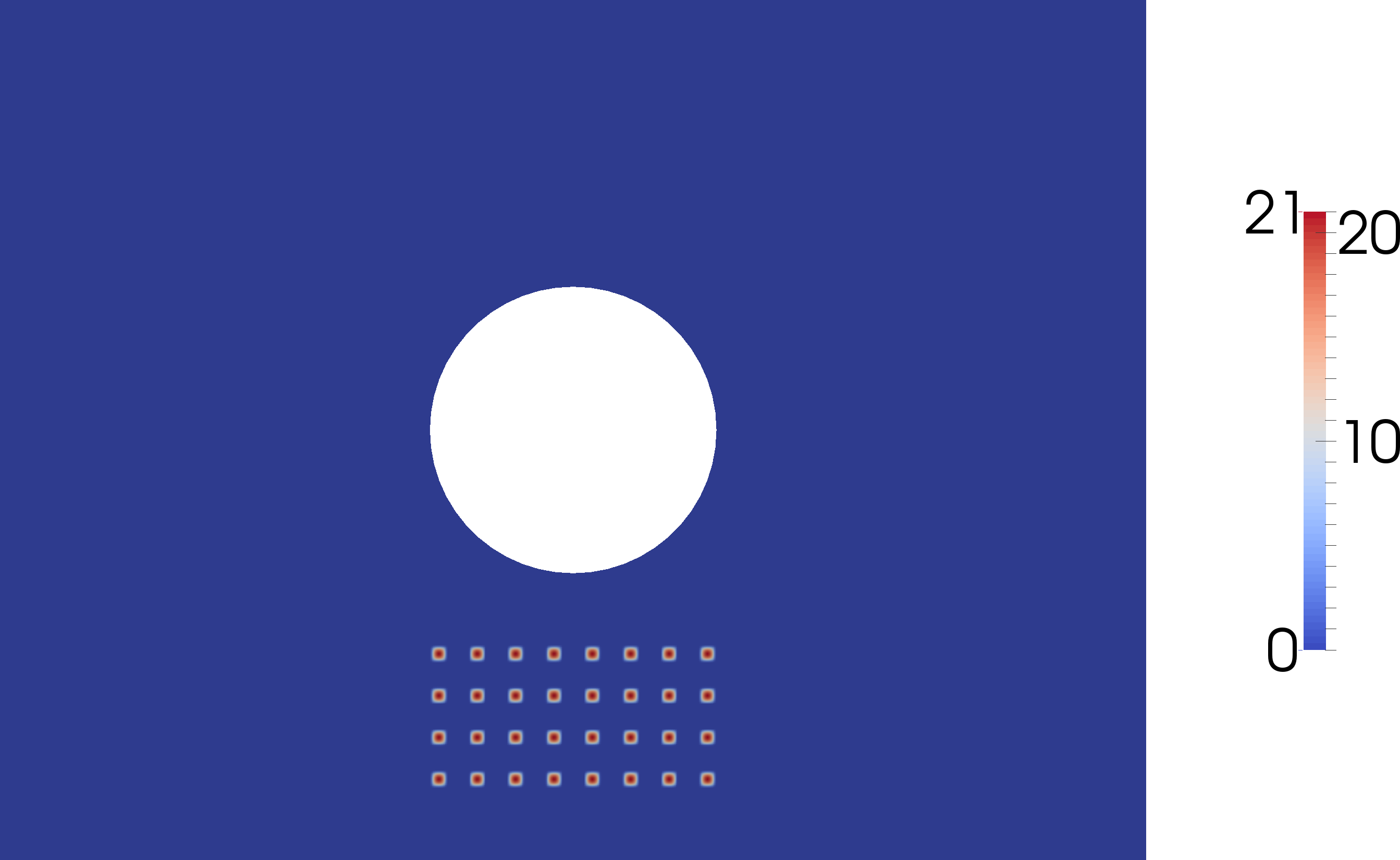}
  \label{fig:results/scenario4_sin/paraview/initial_turbine_friction}
  }
  \subfloat[Optimised turbine positions]{
    \centering
      \includegraphics[width=0.39\textwidth]{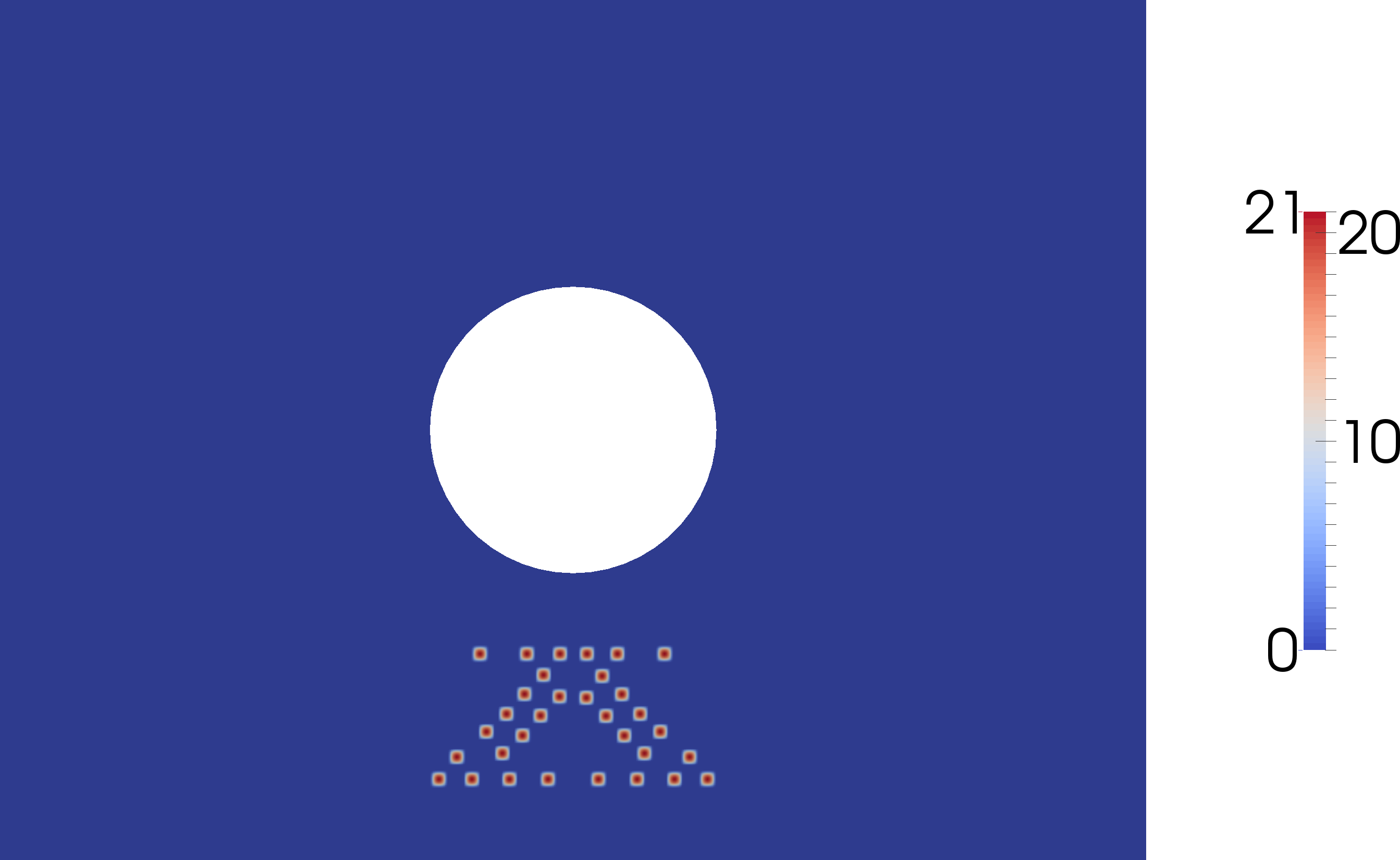}
  \label{fig:results/scenario4_sin/paraview/turbine_friction}
  }
  \\
  \subfloat[Velocity magnitude at the\newline time when the input velocity\newline from the left reaches its maximum]{
    \centering
      \includegraphics[width=0.39\textwidth]{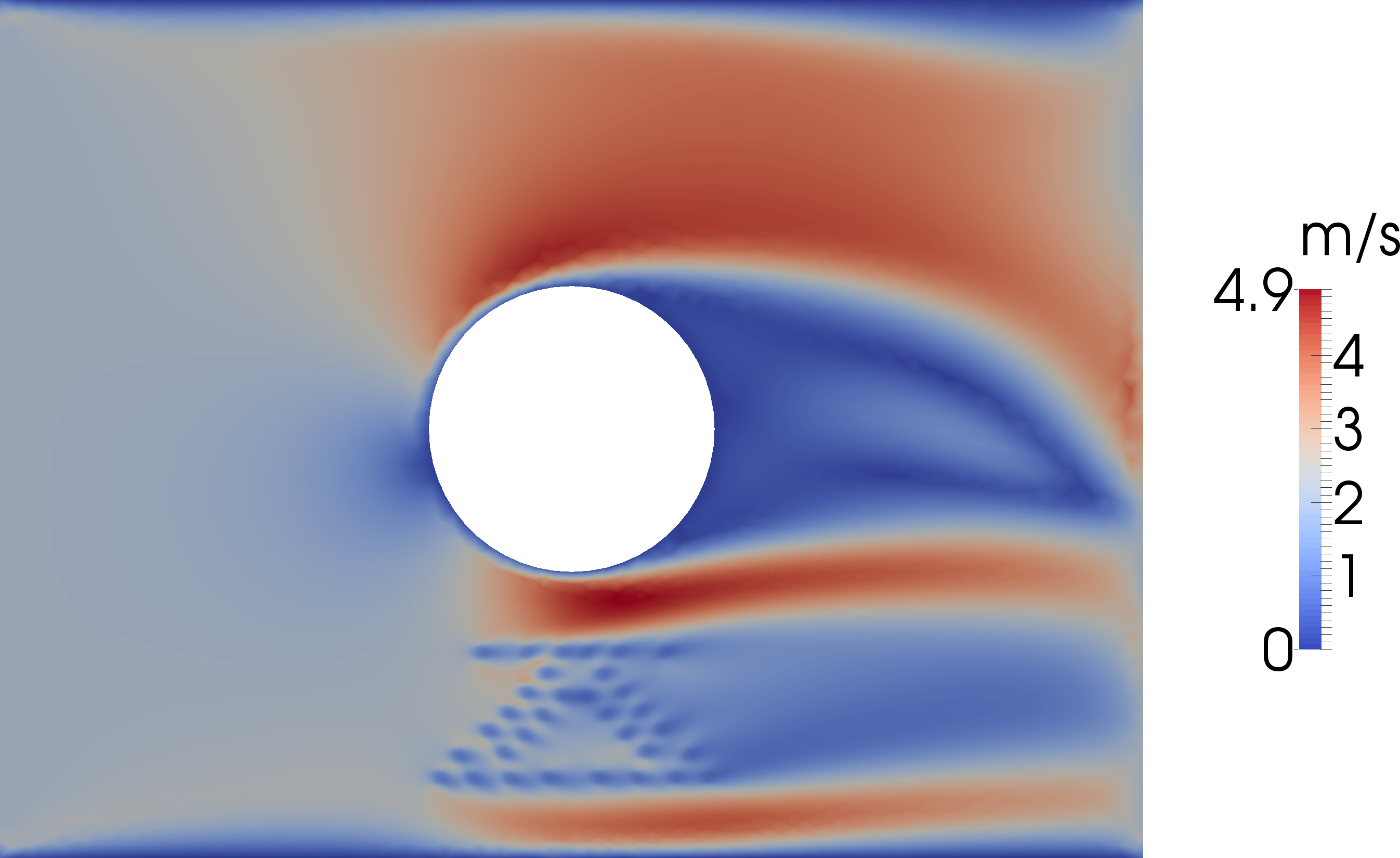}
  \label{fig:results/scenario4_sin/paraview/velocity_left}
  }
  \subfloat[Velocity magnitude at the\newline time when the input velocity\newline from the right reaches its maximum]{
    \centering
      \includegraphics[width=0.39\textwidth]{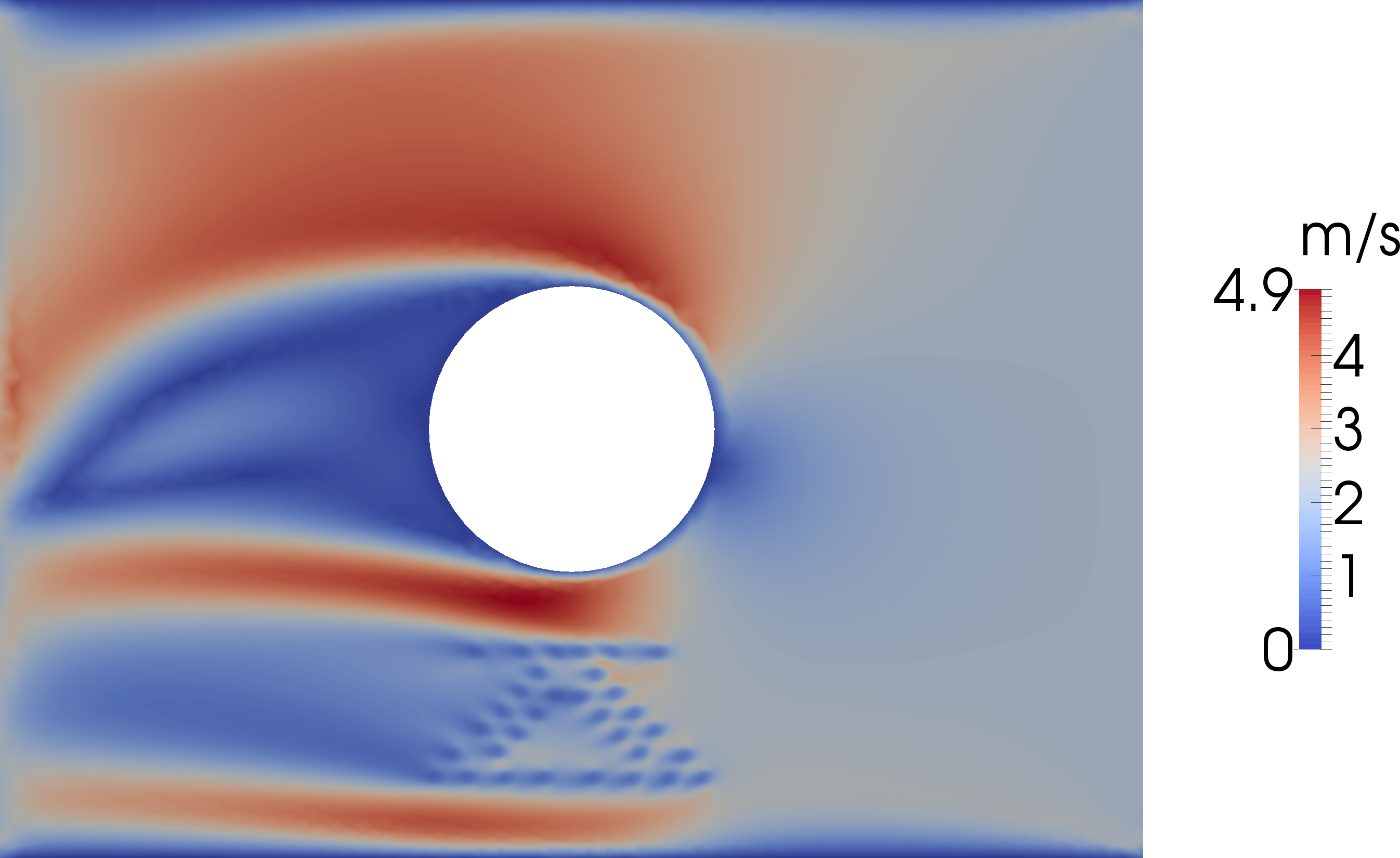}
  \label{fig:results/scenario4_sin/paraview/velocity_right}
  }
  \\
  \subfloat[Optimisation convergence]{
    \centering
      \includegraphics[width=0.39\textwidth]{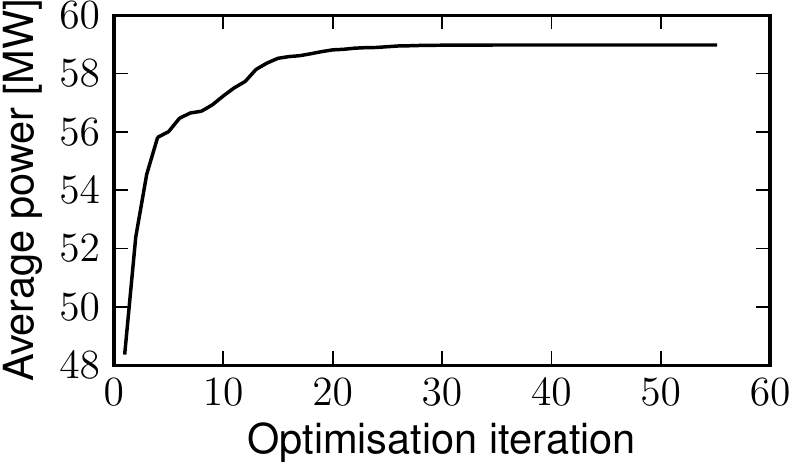}
  \label{fig:results/scenario4_sin/iter_plot}
  }
  \subfloat[Power extraction over time of the optimised configuration]{
    \centering
      \includegraphics[width=0.39\textwidth]{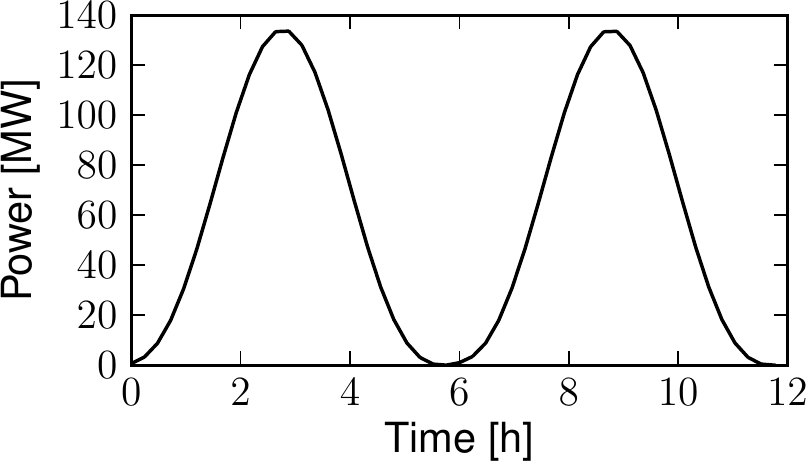}
      \label{fig:results/scenario4_sin/power_timeline} 
      }
  \caption{Results of the non-stationary scenario 4 with minimum distance constraints.}\label{fig:results/scenario4_sin} 
\end{figure}
The final scenario (figure~\ref{fig:scenarios/scenario4}) was solved with the non-stationary shallow water equations.
The simulation time consisted of one $P = 12$~h sinusoidal period with a time step of $\Delta t = 864$~s.
No spin up phase was applied, as its effect is assumed to be small due to the relatively short extent of the domain.
Dirichlet boundary conditions on the left and right boundaries enforced the following sinusoidal in-/outflow velocity:
\begin{equation*}
 u(x, y, t)  =
\begin{pmatrix}
-2\sin\left(2\pi t/P\right) \\
0
\end{pmatrix}.
\end{equation*}
On the remaining boundaries a no-slip boundary condition was applied.

The setup optimised the turbine positions and applied the inequality constraints to enforce a minimum turbine
distance of $30$~m.  The optimisation algorithm terminated after $55$ iterations ($54$ gradient evaluations,
$75$ functional evaluations).  The results are shown in figure~\ref{fig:results/scenario4_sin}. 

The averaged power extracted during one cycle increased by $22\%$ from $48.4$~MW to $59.0$~MW
(figure~\ref{fig:results/scenario4_sin/iter_plot}).  Since the computational domain is symmetric and the
simulation time covered one full period, the optimised layout is expected to be symmetric in the $x$-direction.
The numerical solution, shown in figure~\ref{fig:results/scenario4_sin/paraview/turbine_friction},
indeed shows an almost symmetric result.  The turbine alignment consists of two distorted $\lor$ shapes
whose open ends face the in/-outflow boundaries.  Similar to the previous example, an interpretation of this
alignment is to divert the stream towards the corner of the $\lor$ where turbines can extract large amounts of
power.  An additional row of turbines can be seen parallel to the bottom of the domain.  These turbines are
positioned to capture energy from the flow passing along the boundary.

\section{Farm optimisation in the Inner Sound of the Pentland Firth} \label{sec:pentland}
\begin{figure}[t]
  \centering
  \subfloat[]{
    \centering
    \includegraphics[width=0.39\textwidth]{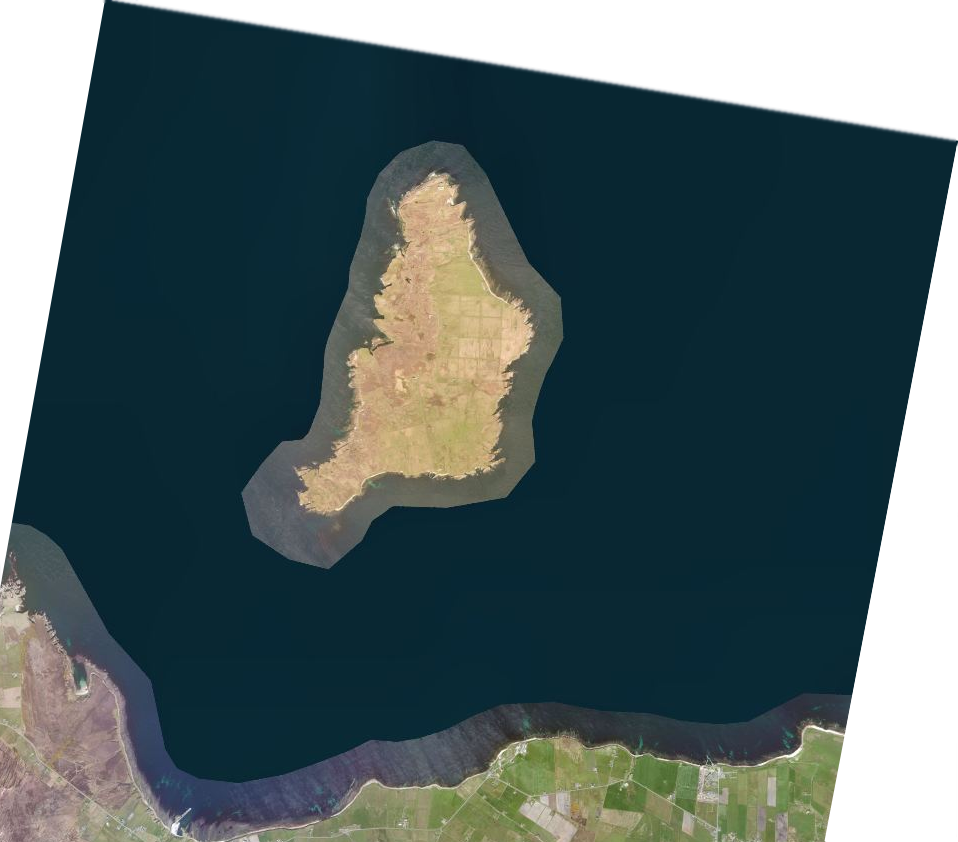}
  \label{fig:pentlandfirth_satellite}
  }
  \subfloat[]{
    \centering
    \includegraphics[width=0.39\textwidth]{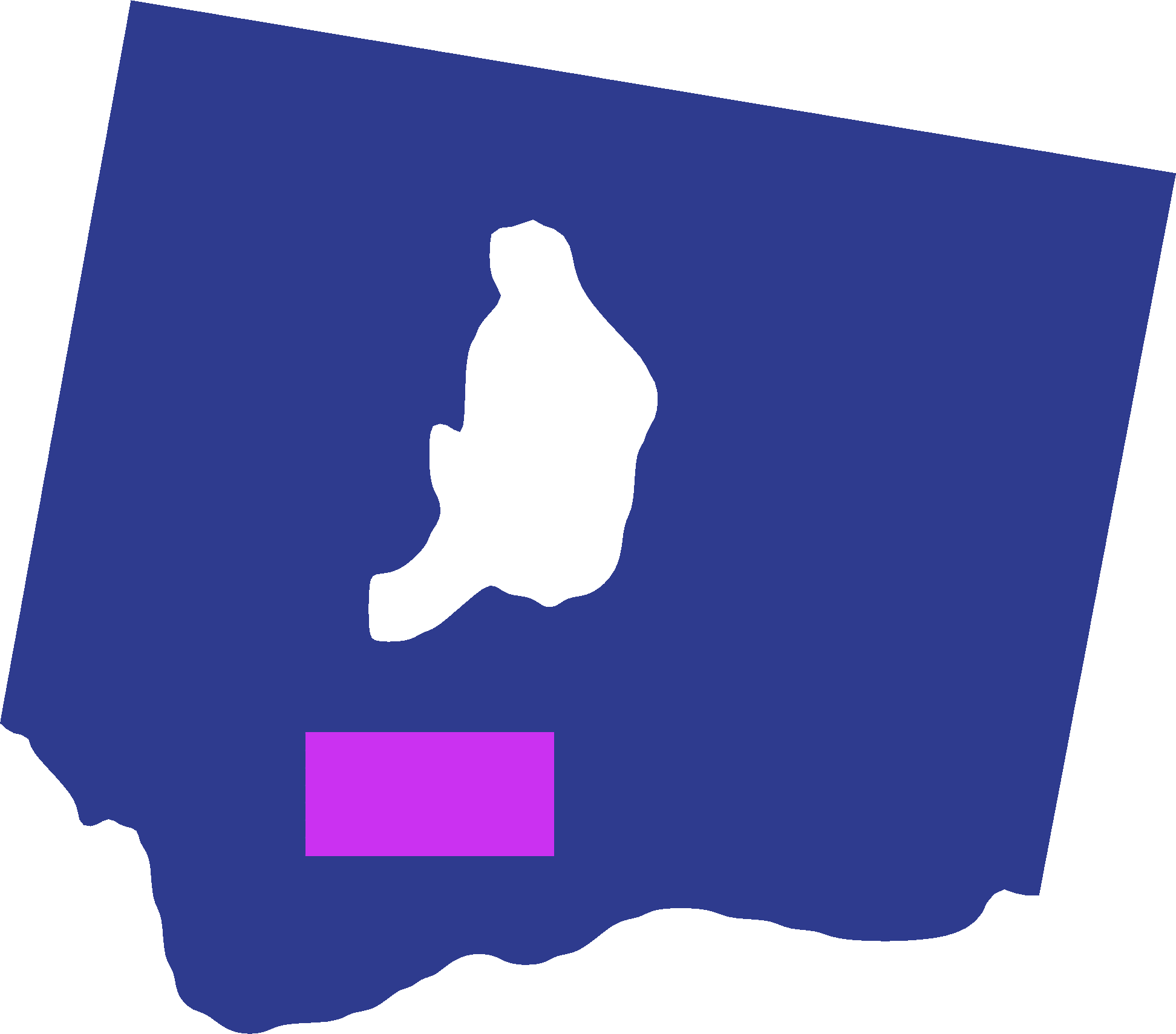}
  \label{fig:pentlandfirth_map}
  }
  \caption{(a): Satellite image of Stroma Island and Caithness (Bing Maps, Microsoft). Satellite imagery is
  only available for the land and nearshore areas. (b): Computational domain with the turbine site marked in pink.}\label{fig:pentlandfirth_domain} 
\end{figure}

A key design feature of the framework is that it should scale to problems in realistic oceanographic domains
with large numbers of turbines in parallel. To demonstrate this capability,
the layout optimisation of a farm in a semi-idealised geometry modelled on the Inner Sound of the Pentland Firth
(figure~\ref{fig:pentlandfirth_satellite}) was conducted on the Stampede supercomputer.  This site is one of
the most promising locations in the UK, and is currently under development by MeyGen Ltd. 

The computational domain is shown in figure~\ref{fig:pentlandfirth_map}.  The pink area marks the turbine site
location which roughly approximates the area used by the MeyGen project.  The discretised domain consists of a
regular mesh in the turbine site area with $2$ m element size, and an unstructured mesh elsewhere with element
sizes ranging from $1.5-200$ m.  The mesh was generated using Gmsh \citep{geuzaine2009} using the GSHHS
shoreline database \citep{wessel1996}. The mesh consists of $1.25 \times 10^6$ elements, which induces a
discretisation with a total of $5.6 \times 10^6$ degrees of freedom with the Taylor-Hood finite element pair.
In this idealised problem, the bathymetry is assumed constant at $H = 50$ m. The
addition of bathymetry is straightforward, and will be presented in future work, but would obscure the
physical interpretation of the results presented here.

\begin{figure}[t]
  \centering
     \subfloat[Optimisation convergence\newline(128 turbines)]
     {
        \centering
          \includegraphics[width=0.39\textwidth]{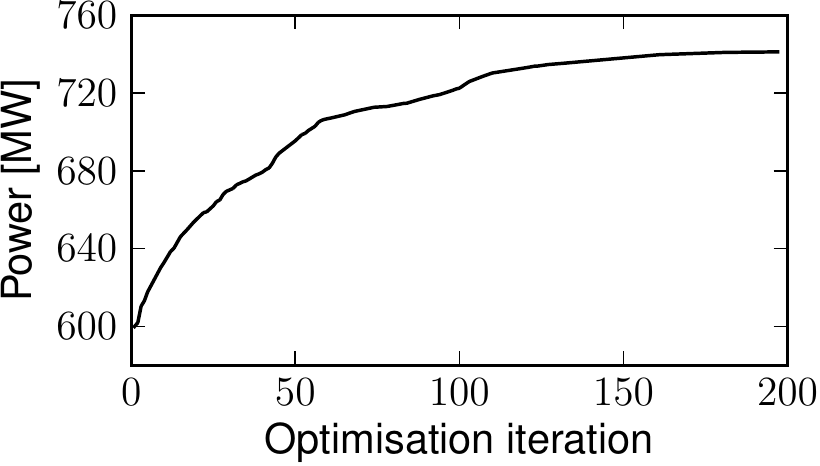}
          \label{fig:pentland_128_iter_plot} 
         }
         \subfloat[Optimisation convergence\newline(256 turbines)]{
        \centering
           \includegraphics[width=0.39\textwidth]{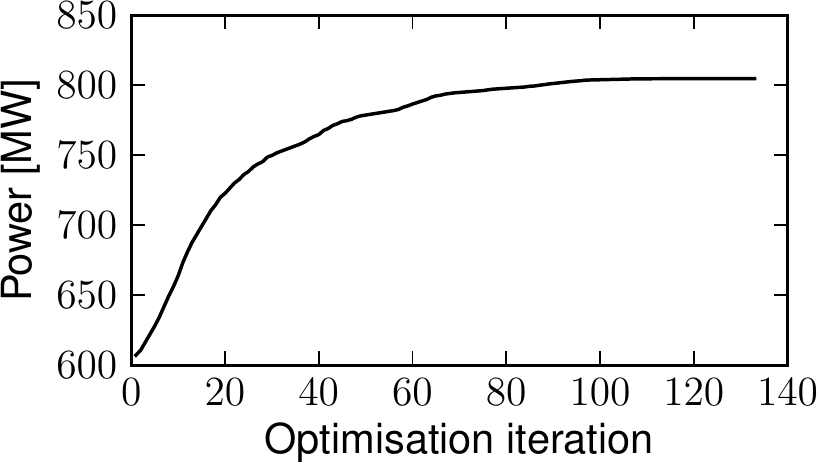}
\label{fig:pentland_256_iter_plot} 
}
   \caption{Convergence of the Pentland Firth optimisation.}
\end{figure}

\begin{figure}[tb]
  \centering
  \subfloat[Initial turbine positions (128 turbines)]{
    \centering
      \includegraphics[width=0.49\textwidth]{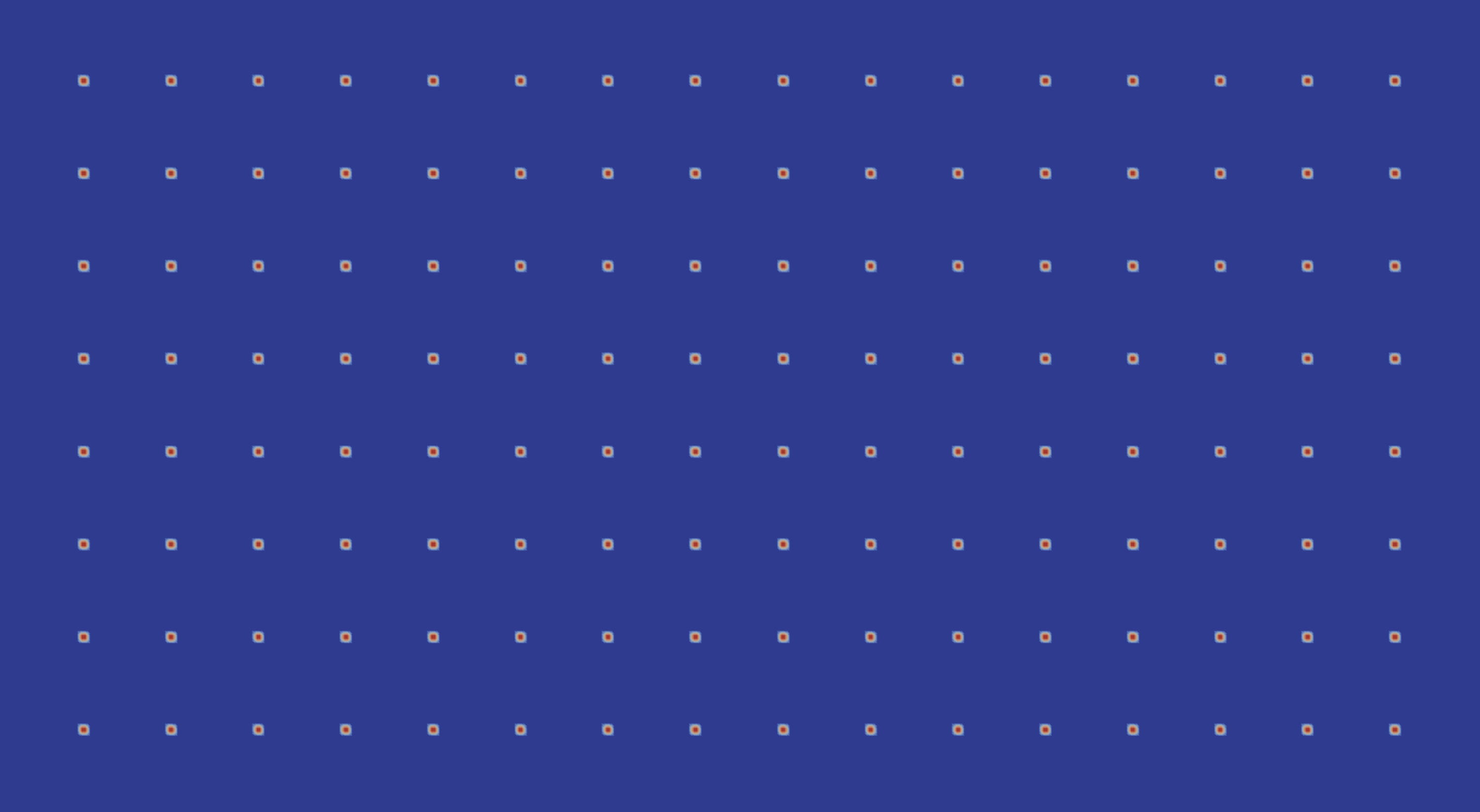}
  \label{fig:pentland_128_initial_position}
  }
  \subfloat[Optimised turbine positions (128 turbines)]{
    \centering
      \includegraphics[width=0.49\textwidth]{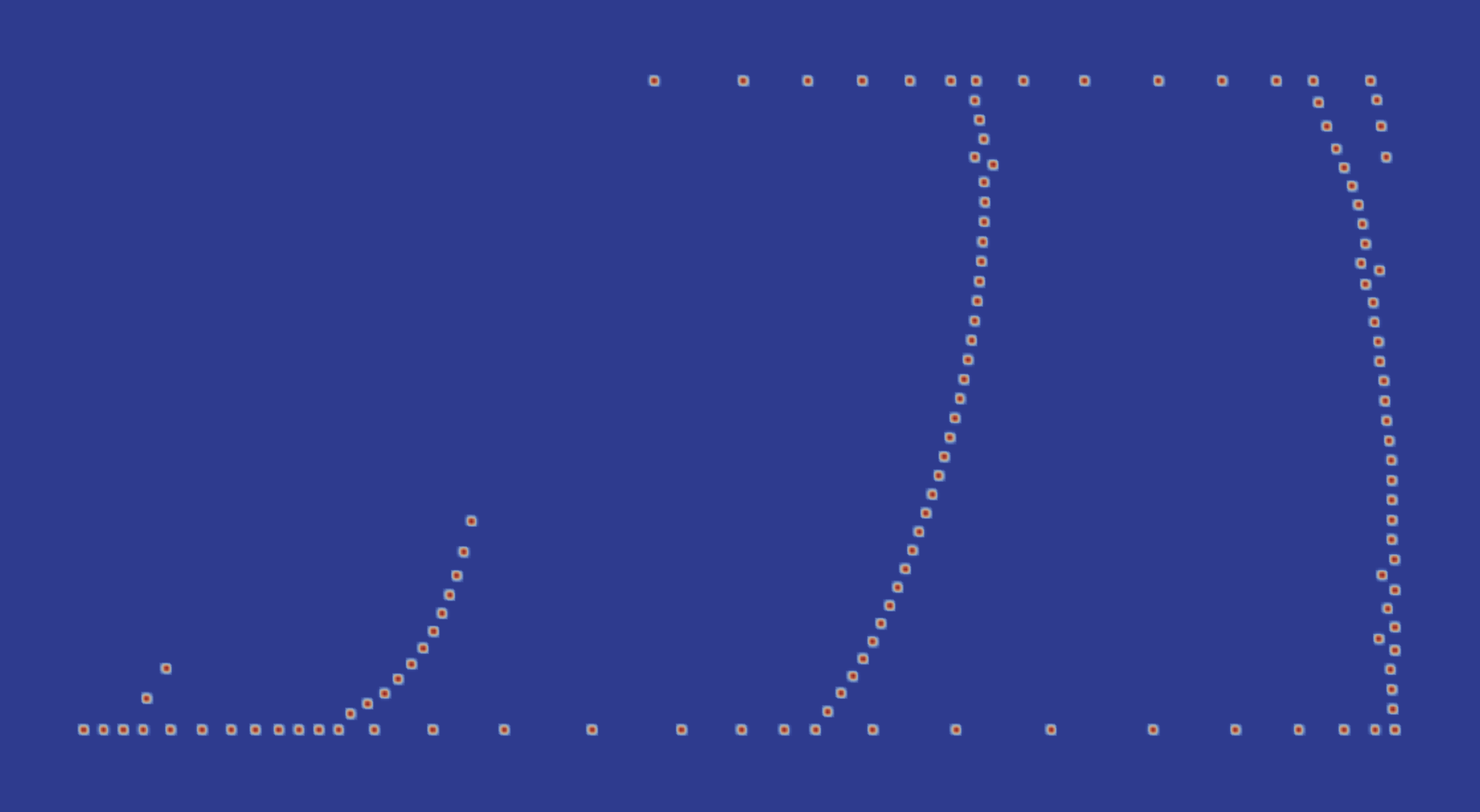}
  \label{fig:pentland_128_final_position}
  }
  \\
  \subfloat[Initial turbine positions (256 turbines)]{
    \centering
      \includegraphics[width=0.49\textwidth]{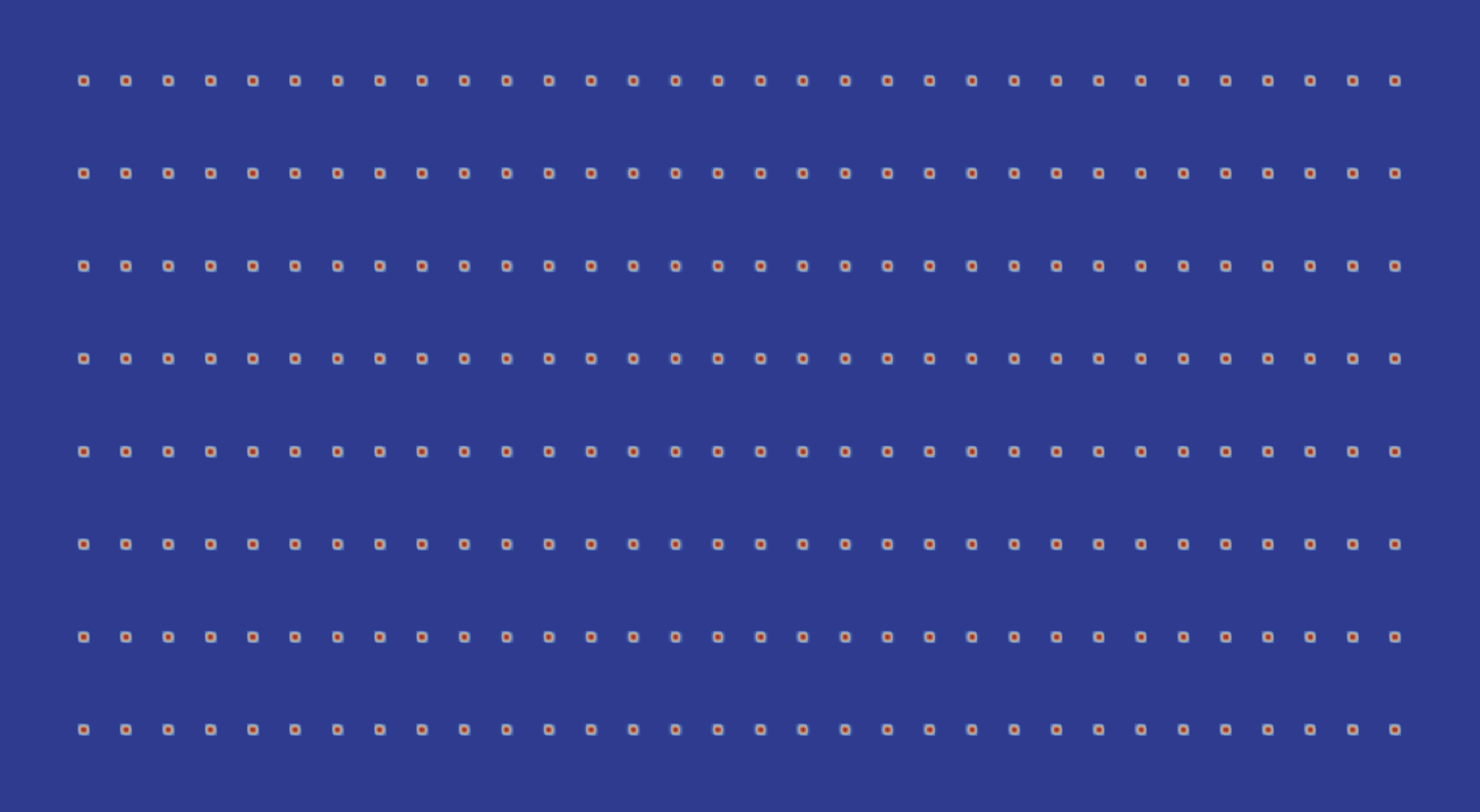}
  \label{fig:pentland_256_initial_position}
  }
  \subfloat[Optimised turbine positions (256 turbines)]{
    \centering
      \includegraphics[width=0.49\textwidth]{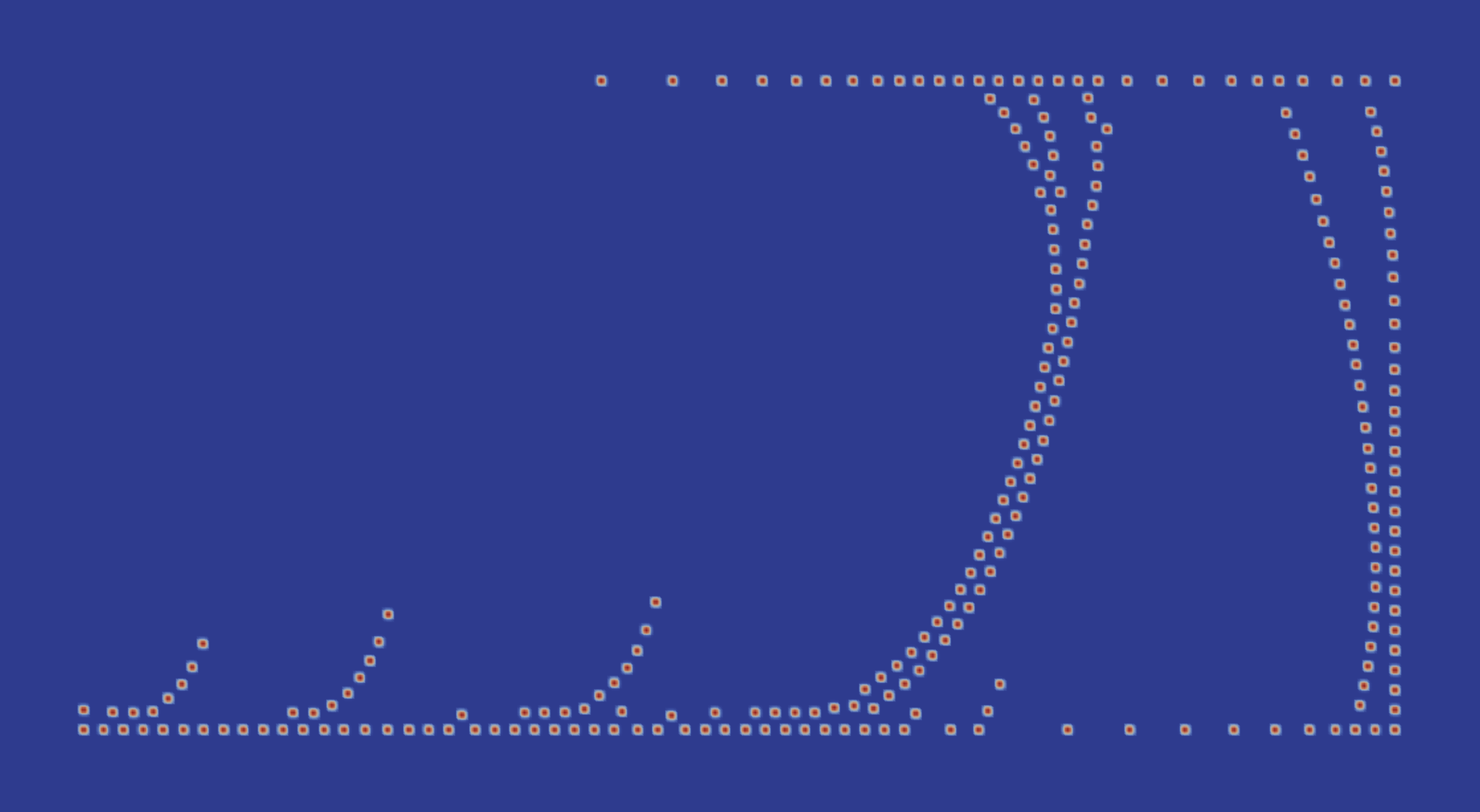}
  \label{fig:pentland_256_final_position}
  }
  \caption{Initial and optimised turbine positions for the Pentland Firth example.}\label{fig:pentlandfirth_turbines}
\end{figure}

\begin{figure}[tb]
  \centering
  \subfloat[Streamline flow visualisation (128 turbines)]{
    \centering
      \includegraphics[width=0.49\textwidth]{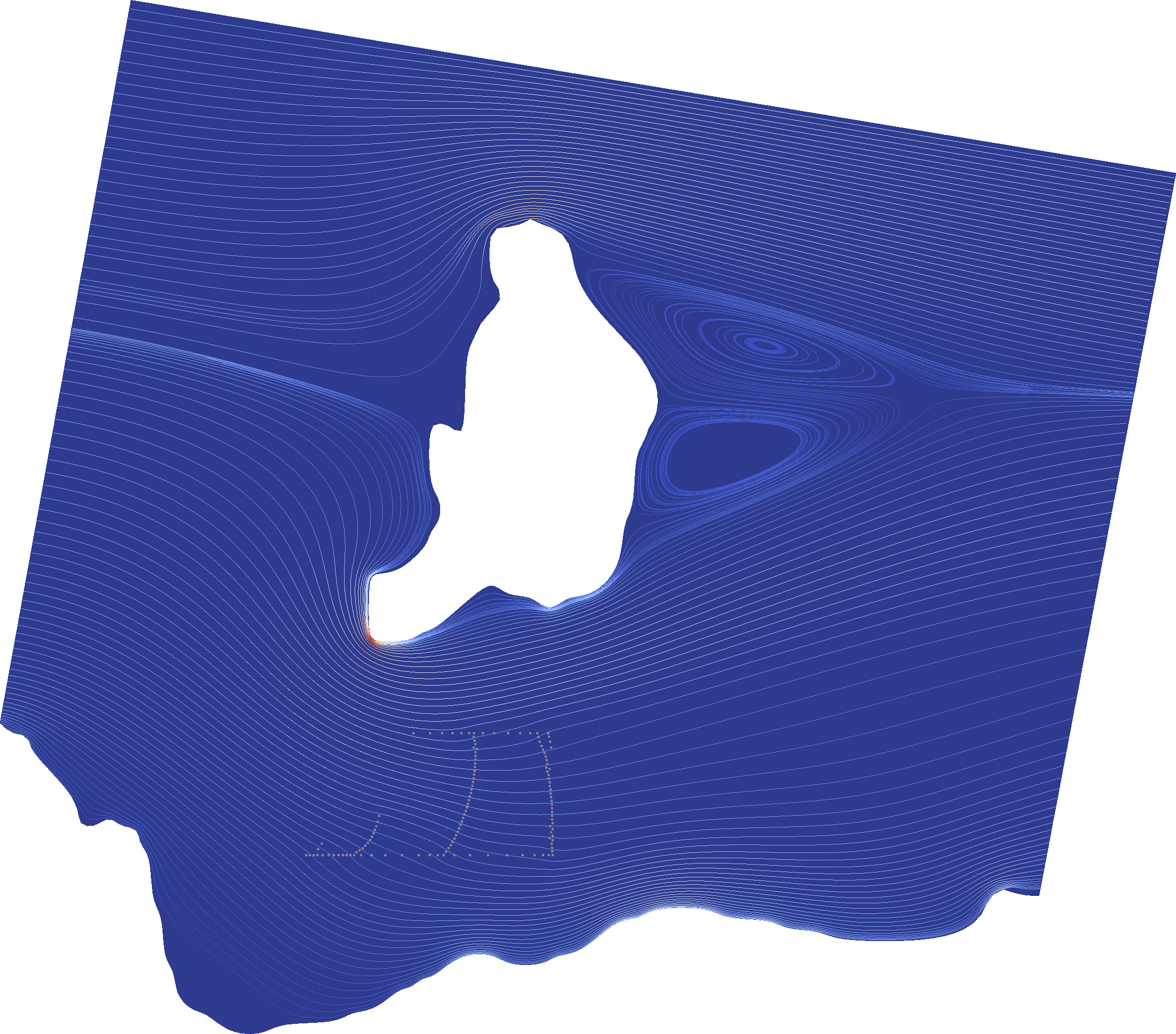}
  \label{fig:pentland_128_streamlines}
  }
  \subfloat[Streamline flow visualisation (256 turbines)]{
    \centering
      \includegraphics[width=0.49\textwidth]{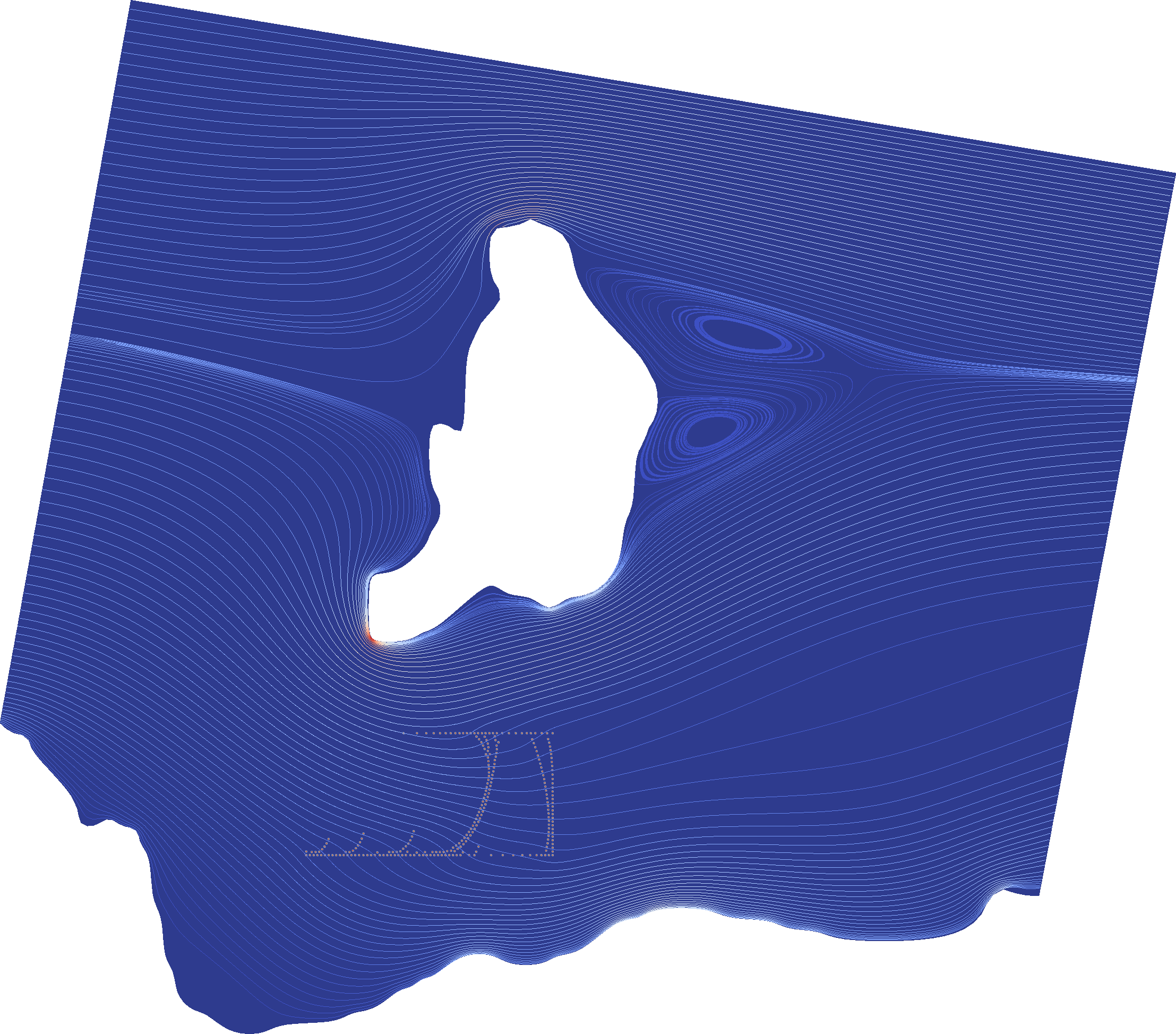}
  \label{fig:pentland_256_streamlines}
  }
  \\
  \subfloat[Turbine power map (128 turbines)]{
    \centering
      \includegraphics[width=0.49\textwidth]{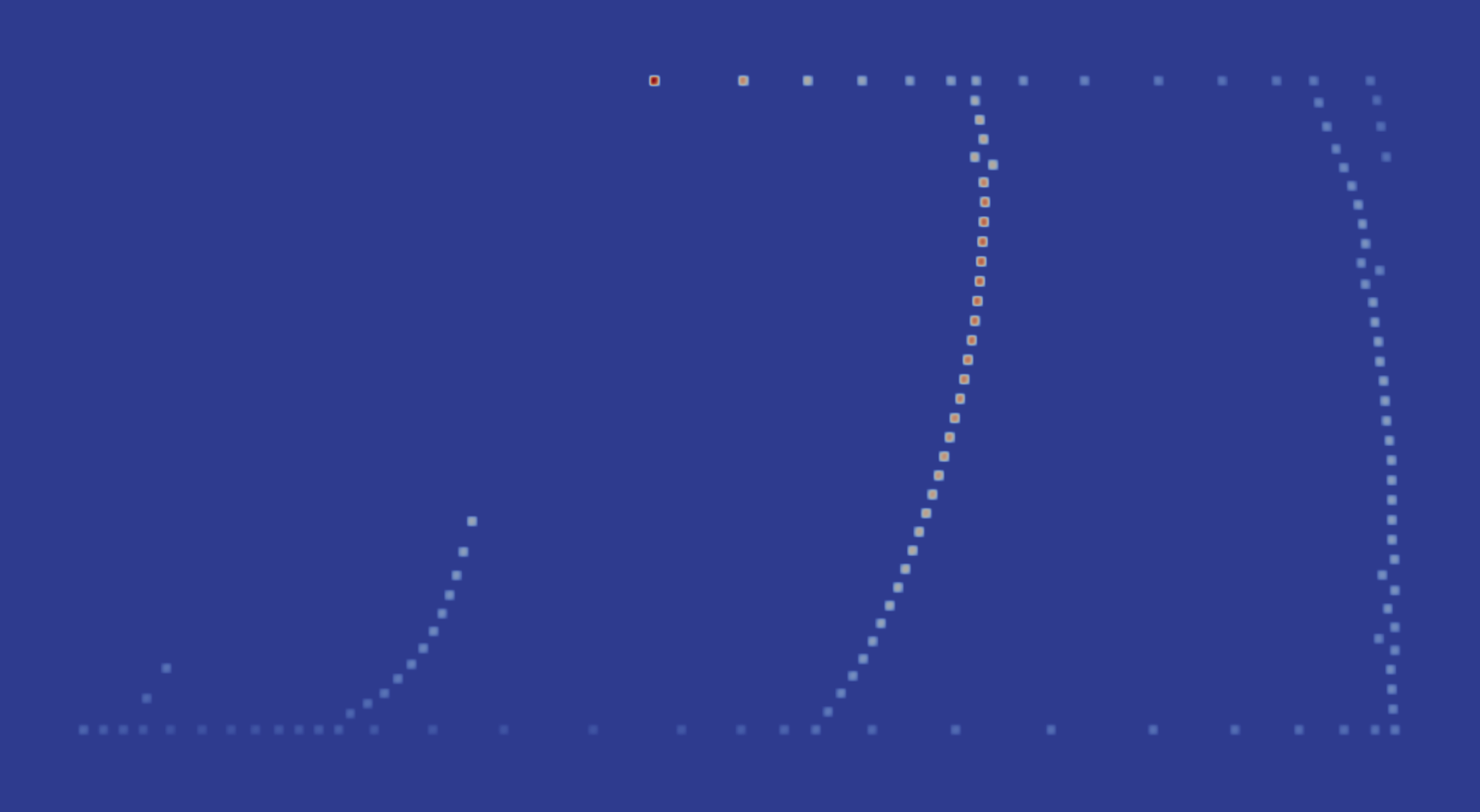}
  \label{fig:pentland_128_power}
  }
  \subfloat[Turbine power map (256 turbines)]{
    \centering
      \includegraphics[width=0.49\textwidth]{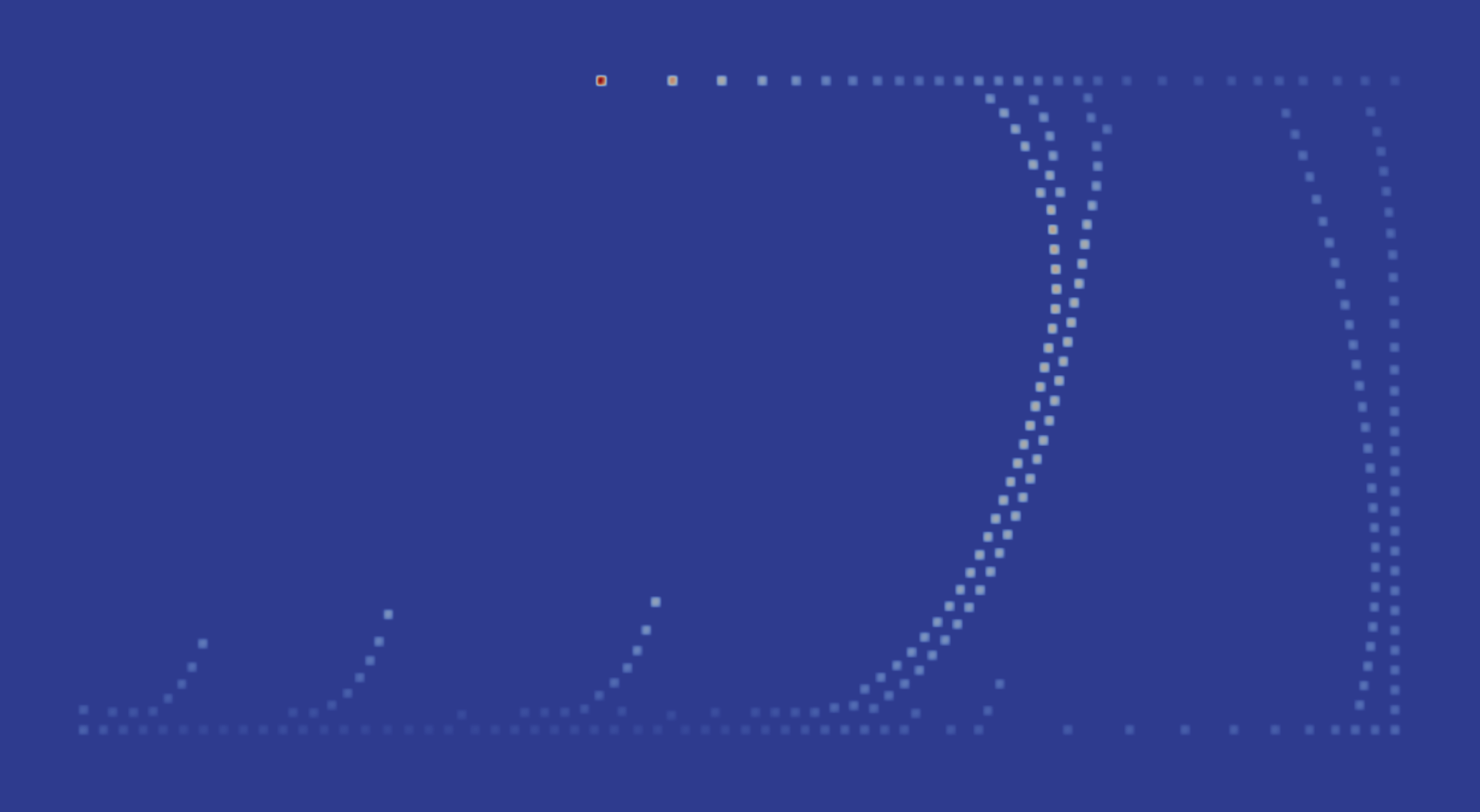}
  \label{fig:pentland_256_power}
  }
  \caption{The streamline and power results for the Pentland Firth example. The power map displays the
  integrand of the functional $J$~(section~\ref{sec:functional}, equation~\eqref{eq:functional_steady}).}
\end{figure}

The farm optimisation was modelled during the flood tide using the stationary shallow water equations.  The
simulations were performed with the same parameter settings as in the previous section
(table~\ref{tab:computational_settings}), but with an increased viscosity coefficient of $\nu = 30\textrm{
m}^2\textrm{s}^{-1}$ and a decreased turbine friction coefficient of $K = 10.5$. This was to ensure the
existence of a steady-state solution. Again, the unsteady case is straightforward and will be presented in
future work. For efficiency reasons, the convergence tolerance of the SLSQP optimisation algorithm was changed
to 1.  The same boundary conditions were used as in the steady-state scenarios of the previous section, with
the inflow condition imposed on the western boundary and the outflow condition enforced on the eastern
boundary. 

Two optimisation runs were conducted, with 128 and 256 turbines respectively. Both were performed on
the Stampede supercomputer at the Texas Advanced Computing Center on 64 cores; both took between 24
and 48 hours of real time to complete (one functional evaluation took approximately 270 seconds, one gradient evaluation took
approximately 90 seconds). 
The 128 turbine case converged
in 197 iterations (197 gradient evaluations, 349 functional evaluations), and increased the power extraction
from 600 MW to 741 MW, an increase of 24\% (figure \ref{fig:pentland_128_iter_plot}).  The 256 turbine case
converged in 133 iterations (133 gradient evaluations, 185 functional evaluations), and increased the power
extraction from 607 MW to 804 MW, an increase of 33\% (figure \ref{fig:pentland_256_iter_plot}). 
Even though eight times as many turbines are to be optimised compared to the previous section, the number of iterations
has hardly increased at all, suggesting the suitability of the method for larger arrays.

The initial and optimised configurations are shown in figure~\ref{fig:pentlandfirth_turbines}.  In both cases,
the optimisation algorithm builds very similar structures.  The key features of the optimised layouts are the following: 
\begin{itemize}
 \item Walls of turbines aligned with the north and south boundaries; the southern wall extends for the entire
 zonal length, while the northern wall only extends for the eastern two-thirds of the site.  The turbine density
 on the western side of the northern wall appears to decay in a similar way in both cases.
 \item Eastern and western `barrages' of densely packed turbines; in the 128 turbine case, the barrages consist
 of a single meridional column, while for the 256 turbine case the extra turbines are used to form additional
 columns. The optimisation algorithm automatically staggers these extra columns to minimize wake shadowing
 (figure~\ref{fig:pentland_256_final_position}).  These
 barrages are aligned perpendicular to the flow field (figures~\ref{fig:pentland_128_streamlines}
 and~\ref{fig:pentland_256_streamlines}) and extract the majority of the power (figures~\ref{fig:pentland_128_power} and~\ref{fig:pentland_256_power}).
 \item `Spurs' arcing from the southern wall in a north-easterly direction; in the 128 turbine case two spurs
 can be seen, while there are three spurs in the 256 turbine case.  Again, the optimisation aligns the spurs
 perpendicular to the flow field. The spur length is chosen such that the majority of streamlines only
 intersect with two rows of turbines (figures~\ref{fig:pentland_128_streamlines}
  and~\ref{fig:pentland_256_streamlines}).
\end{itemize}

We hypothesise that these structures serve the following functions:
\begin{itemize}
 \item The main purpose of the southern and northern walls is to funnel the flow into and retain the flow
   inside the site boundary. Similar features are found in all scenarios of the previous section. The water
   predominantly flows in from the northwest, which is why the northern wall stops short on the western side.
   We hypothesise that the decay of the turbine density on the western side of the northern wall acts to
   funnel water through the barrages. A similar density decay is visible
   in the optimised layout of scenario 1 (figure~\ref{fig:results/scenario1_no_ineq/paraview/turbine_friction}).
 \item The flow trapped inside the site domain attempts to escape through the northern and southern walls,
   which causes the arcing of the western and eastern barrages close to the northern and southern boundaries
   (figures~\ref{fig:pentland_128_streamlines} and~\ref{fig:pentland_256_streamlines}).
   Since the prevailing incoming flow is towards the southeast, more turbines are placed on the southern wall to retain
   the flow. This motivates the deployment of the spurs on the western side of the southern wall. 
\end{itemize}

It is not physically meaningful to compare the optimised power extractions for the 128 and 256 turbine cases, as a realistic 
power curve was not used. The maximum velocity in the site for the 128 turbine case was $3.7$~ms$^{-1}$, while it was 3.0~ms$^{-1}$ for the 256 turbine case. 
Due to the cubic dependence of the power extraction on the speed, the power per turbine is approximately doubled in the 128 turbine case, and
so the total power extraction of the two cases are almost the same. However, if the turbine model enforced a rated speed beyond which
no extra power was extracted, this approximate equality would not hold.

\section{Conclusions}
In this work, the optimal configuration of tidal turbine farms was formulated as an optimisation problem constrained
by partial differential equations describing the flow. This formulation allows the direct application of sophisticated mathematical
techniques to its solution, particularly the adjoint technique for rapidly evaluating gradients. The optimisation of tidal turbine farms is necessary to realise
their full potential, but without gradient-based optimisation techniques the computational cost is prohibitive \citep{vennell2012c}.

The approach presented here has several key advantages. It fully accounts for the nonlinear interactions
between the geometry, the turbines, and the flow throughout the optimisation. The use of gradient-based
optimisation algorithms combined with the adjoint technique enables the use of physically-realistic flow
models, even for a large number of turbines. Once the flow model inputs are specified (domain, boundary
conditions, initial array configuration, etc.), the optimisation is fully automatic. The approach extends naturally to more realistic flow models,
and to different functionals of interest such as profit or environmental impact.

The algorithm was first applied to the optimisation of four idealised scenarios, both to demonstrate the
capability of the method and to build physical intuition. In all cases, the optimisation algorithm was
successful in significantly increasing the power extracted by the farm, at a computationally feasible cost.
The algorithm was then successfully applied to a more realistic optimisation problem, involving a site of
major industrial interest, accurate shoreline geometry, and an industrially relevant number of turbines.

Four main extensions are required to apply this in an industrial setting. Firstly, the simulations
must be driven by realistic tidal forcing, and incorporate a wider model domain around the site to be
investigated. 
This also includes extending the forward model to be forced by a head loss instead of a fixed inflow velocity 
to incorporate the impact of large arrays on the free-stream velocity~\citep{vennell2010}. 
Secondly, bathymetric effects must be accounted for. Thirdly, the
flow model must then be validated against real-world measurements. Finally, the wake modelling should
be improved via a turbulence closure and a realistic power curve used. All of these advances are the
subject of ongoing work.

The source-code of the turbine farm optimisation software and all examples are 
open-source and available at \url{http://opentidalfarm.org}.

\section{Acknowledgements}
This work is supported by the Grantham Institute for Climate Change, a Fujitsu CASE studentship, EPSRC grants EP/I00405X/1, EP/J010065/1, EP/K503733/1 and EP/K030930/1, an EPSRC Pathways to Impact award, and a Center of Excellence grant from the Research Council of Norway to the Center for Biomedical Computing at Simula Research Laboratory. High-performance Computing support was provided by the Texas Advanced Computing Center. The authors would like to
acknowledge helpful discussions with S. C. Kramer, A. S. Candy, A. Avdis of Imperial College London, and R. Caljouw and S. Crammond of MeyGen Ltd.
The authors would also like to thank R. Vennell of the University of Otago for his extremely thorough and constructive review. 

\bibliographystyle{elsarticle-harv}
\bibliography{literature}
\end{document}